# Some Properties of Mordukhovich Derivatives with Applications to Locally Variational Inequalities in Banach Spaces

Jinlu Li

Department of Mathematics
Shawnee State University
Portsmouth, Ohio 45662 USA
jli@shawnee.edu

**Abstract**. In this paper, we will investigate the connection between Fréchet derivatives and Mordukhovich derivatives of single-valued mappings in Banach spaces. We will find some properties of Mordukhovich derivatives of set-valued mappings, which will be demonstrated by the set-valued metric projection operator in finite dimensional Banach spaces. We introduce the concept of locally variational inequalities in Banach spaces, which focus on finding solutions within some specific neighborhoods rather than across the entire global domain; and we find the connection between the solutions of locally variational inequality problems and the Mordukhovich derivatives.

## 1. Introduction and Preliminaries

In Banach analysis with respect to single-valued mappings, the most popular concepts of differentiation are Gâteaux directionally differentiability and Fréchet differentiability (See $[5-9, 13, 26-28]$). In particular, in $[10-12, 17, 19]$, the explicit Gâteaux and Fréchet derivatives of some special mappings, such as the metric projection operator, the normalized duality mapping, have been calculated, which are used for calculate their Mordukhovich derivatives.

In modern analysis in Banach spaces, the theory of Mordukhovich differentiation for set-valued mappings in Banach spaces forms the foundation of generalized differentiation in set-valued and variational analysis in Banach spaces (See $[22-25]$). The theory of Mordukhovich differentiation has been widely applied to both pure and applied mathematics such as optimization theory, game theory, approximation theory, and so forth (See [1−4, 16, 21−25]).

When single-valued mappings are considered as special cases of set-valued mappings with values being singletons, the Mordukhovich derivatives of the metric projection operator and the normalized duality mapping have been calculated by using the connection between Fréchet and Mordukhovich derivatives, which is given in Theorem 1.38 in [23] (See [10−12, 14−16, $18, 20-25$]).

Before we review the concepts of Mordukhovich derivatives for set-valued mappings, we first review the concepts and some properties of Fréchet derivatives for single-valued mappings in Banach spaces.

Let $(X, \|\cdot\|_X)$ and $(Y, \|\cdot\|_Y)$ be real Banach spaces with topological dual spaces $(X^*, \|\cdot\|_{X^*})$ and $(Y^*, \|\cdot\|_{Y^*})$, respectively. Let $\langle\cdot, \cdot\rangle_X$ denote the real canonical pairing between $X^*$ and $X$ and $\langle\cdot, \cdot\rangle_Y$ the real canonical pairing between $Y^*$ and $Y$. Let $\theta_X$ and $\theta_Y$ denote the null elements of $X$ and $Y$, respectively. Let $\theta_{X^*}$ and $\theta_{Y^*}$ denote the null elements of $X^*$ and $Y^*$, respectively. Let $A$ be a nonempty convex and open subset of $X$. Let $f: A \to Y$ be a single-valued mapping. Let $\bar{x} \in A$. If there is a continuous and linear mapping

$\nabla f(\bar{x})$: $X \to Y$ such that,

$$\lim_{u \to \theta_X} \frac{f(\bar{x}+u)-f(\bar{x})-\nabla f(\bar{x})(u)}{\|u\|_X} = \theta_Y, \tag{1.1}$$

then $f$ is said to be Fréchet differentiable at $\bar{x}$ and $\nabla f(\bar{x})$ is called the Fréchet derivative of $f$ at $\bar{x}$.

Then, we review the concepts and some properties of Mordukhovich derivatives for set-valued mappings in Banach spaces (See [21 −26]).

Let $F$: $A \rightrightarrows Y$ be a set valued mapping. The graph of $F$ is defined by the following subset in $A \times Y$

$$\text{gph}F = \{(x, y) \in A \times Y : y \in F(x)\}.$$

For $(x, y) \in \text{gph}F$, that is, for $x \in A$ and $y \in F(x)$, the Mordukhovich derivative (which is also called Fréchet coderivative, Mordukhovich coderivative, or just coderivative) of $F$ at point $(x, y)$ is a set-valued mapping $\widehat{D}^*F(x, y)$: $Y^* \rightrightarrows X^*$. By Definitions 1.13 and 1.32 in Chapter 1 in [23], $\widehat{D}^*F(x, y)$ is defined, for any $y^* \in Y^*$, by

$$\widehat{D}^*F(x, y)(y^*) = \left\{x^* \in X^*: \limsup_{\substack{(u,v)\to(x,y)\\ u\in A \text{ and } v\in F(u)}} \frac{\langle x^*, u-x\rangle_X - \langle y^*, v-y\rangle_Y}{\|u-x\|_X + \|v-y\|_Y} \leq 0\right\}$$

$$= \left\{x^* \in X^*: \limsup_{\substack{(u,v)\to(x,y)\\ (u,v) \in \text{gph}F}} \frac{\langle x^*, u-x\rangle_X - \langle y^*, v-y\rangle_Y}{\|u-x\|_X + \|v-y\|_Y} \leq 0\right\}. \tag{1.2}$$

For $(x, y) \in \text{gph}F$, if there is $y^* \in Y^*$ such that

$$\widehat{D}^*F(x, y)(y^*) \neq \emptyset, \text{ for some } y^* \in Y^*, \tag{1.3}$$

then, $F$ is said to be Mordukhovich differentiable at $(x, y)$. If $(x, y) \notin \text{gph}F$, then, we define

$$\widehat{D}^*F(x, y)(y^*) = \emptyset, \text{ for any } y^* \in Y^*.$$

By the above definition (1.2), $\widehat{D}^*F(x, y)$: $Y^* \rightrightarrows X^*$ is a set-valued mapping and it is called the Mordukhovich derivative (or the Mordukhovich coderivative, or coderivative) of $F$ at $(x, y)$.

In particular, let $f$:$A \to Y$ be a single valued mapping. Let $x \in X$ and $y \in Y$ with $y = f(x)$, the Mordukhovich derivative (coderivative) of $f$ at point $(x, y)$ is also a set-valued mapping, denoted by $\widehat{D}^*f(x, y)$: $Y^* \rightrightarrows X^*$. The definition (1.2) becomes that,

$$\widehat{D}^*f(x, y)(y^*) = \left\{z^* \in X^*: \limsup_{\substack{(u,f(u))\to(x,y)\\ u\in A}} \frac{\langle z^*, u-x\rangle_X - \langle y^*, f(u)-y\rangle_Y}{\|u-x\|_X + \|f(u)-y\|_Y} \leq 0\right\}, \text{ for any } y^* \in Y^*.$$

More strictly, if $f$:$A \to Y$ is a continuous single-valued mapping, then

$$\widehat{D}^*f(x,y)(y^*) = \left\{z^* \in X^*: \limsup_{\substack{u\to x\\ u\in A}} \frac{\langle z^*, u-x\rangle_X - \langle y^*, f(u)-y\rangle_Y}{\|u-x\|_X + \|f(u)-y\|_Y} \leq 0\right\}, \text{ for any } y^* \in Y^*. \tag{1.4}$$

The following theorem shows the connection between Fréchet derivatives and Mordukhovich derivatives for single-valued mappings, which provides a powerful tool to calculate the Mordukhovich derivatives by the Fréchet derivatives of single-valued mappings between Banach spaces.

**Theorem 1.38 in [23]**. *Let X be a Banach space with dual space $X^*$and let $f$:X ⟶ Y be a single-valued mapping. Suppose that $f$ is Fréchet differentiable at $x \in X$ with $y = f(x)$. Then, the Mordukhovich derivative of $f$ at $x$ satisfies the following equation*

$$\widehat{D}^*f(x,y)(y^*) = \{(\nabla f(x))^*(y^*)\}, \ \textit{for all } y^* \in Y^*. \tag{1.5}$$

Under the sense of Mordukhovich differentiability of *F* given by (1.3), for single-valued mapping *f*, the results (1.5) of Theorem 1.38 in [23] roughly speaks that, for any $(x, y) \in \text{gph}f$,

$$f \text{ is Fréchet differentiable at } x \implies f \text{ is Mordukhovich differentiable at } (x,y). \tag{1.6}$$

Furthermore, to investigate more details of the connection between Mordukhovich derivatives and Fréchet derivatives for single-valued mappings between Banach spaces, in section 2, we provide some counter examples to show that the converse of (1.6) does not hold. That is, we will prove that, there are some points $(x, y) \in \text{gph}f$ such that

$$f \text{ is Mordukhovich differentiable at } (x,y) \not\Rightarrow f \text{ is Fréchet differentiable at } x. \tag{1.7}$$

In section 3, we prove some properties of Mordukhovich derivatives of set-valued mappings. In section 4, in particularly, we investigate the closeness and convexity of the Mordukhovich derivatives of set-valued mappings, which are also demonstrated by some examples.

In section 5, we introduce the concept of locally variational inequality problems in Banach spaces, which focus on finding solutions within some specific neighborhoods of the considered variables rather than across the entire global domain. We find the connection between solutions of locally variational inequality problems and Mordukhovich derivatives and provide some applications of Mordukhovich derivatives to solving locally variational inequality problems.

## 2. Connection between Fréchet and Mordukhovich Differentiability for Single-valued Mappings

Throughout this paper, we always let ℕ, ℝ and $\mathbb{R}_+$ respectively denote the set of nonnegative integers, the set of real numbers and the set of nonnegative real numbers. For any positive integer *n*, let $\mathbb{R}^n$ denote the *n*-d Euclidean space. In particular, $\mathbb{R}^1 = \mathbb{R}$. Let 𝕊 denote the linear space of real sequences with term-by-term addition and scalar multiplication. For *n* > 1, let $\|\cdot\|_n$ denote the Hilbert norm on $\mathbb{R}^n$.

In this section, in order to show (1.7) under the sense (1.3), more precisely, we will construct some counter examples to prove that the converse of (1.6) is not true under the following sense: for a single-valued mapping $f$:*X* ⟶ *Y*, there are some points $(x, y) \in \text{gph}f$ and $y^* \in Y^*$ satisfying

$$\widehat{D}^*f(x,y)(y^*) \neq \emptyset \quad \text{and} \quad \nabla f(x) \text{ does not exist.}$$

For this purpose, we start with a simple Banach space ℝ.

**Example 2.1**. Let $X = Y = X^* = Y^* = \mathbb{R}$. Define $f: \mathbb{R} \to \mathbb{R}_+$ by

$$f(x) = |x|, \text{ for any } x \in \mathbb{R}.$$

Then, by $f(0) = 0$, the Fréchet differentiability of $f$ at point 0, and the Mordukhovich differentiability of $f$ at point (0,0) have the following properties.

(a) $f$ is not Fréchet differentiable at 0, that is,

$$\nabla f(0) \text{ does not exist.} \tag{2.1}$$

(b) For any $y^* \geq 0$,

$$\widehat{D}^* f(0,0)(y^*) = [-y^*, y^*]. \tag{2.2}$$

(c) For any $y^* < 0$,

$$\widehat{D}^* f(0,0)(y^*) = \emptyset.$$

*Proof*. Proof of (a). Assume that $\nabla f(0)$ exists, which is a continuous and linear operator $\nabla f(0): \mathbb{R} \to \mathbb{R}$. Hence, $\nabla f(0)$ is represented by a real number. By (1.3), we have

$$\begin{aligned} 0 &= \lim_{u \to 0} \frac{f(0+u) - f(0) - \nabla f(0)(u)}{|u|} \\ &= \lim_{u \to 0} \frac{|u| - \nabla f(0)(u)}{|u|} \\ &= \lim_{u \to 0} \left(1 - \nabla f(0) \frac{u}{|u|}\right) \\ &= \begin{cases} 1 - \nabla f(0), & \text{if } u \downarrow 0, \\ 1 + \nabla f(0), & \text{if } u \uparrow 0. \end{cases} \end{aligned}$$

This contradiction proves (2.1) for part (a).

Proof of (b). Let $y^* \geq 0$. Since $f$ is a continuous single-valued mapping, for any $x^* \in [-y^*, y^*]$, by (1.4), we calculate

$$\begin{aligned} &\limsup_{u \to 0} \frac{\langle x^*, u - 0\rangle_X - \langle y^*, f(u) - 0\rangle_Y}{\|u - 0\|_X + \|f(u) - 0\|_Y} \\ &= \limsup_{u \to 0} \frac{\langle x^*, u\rangle_X - \langle y^*, |u|\rangle_Y}{|u| + |u|} \\ &= \limsup_{u \to 0} \frac{x^* u - y^*|u|}{2|u|} \\ &\leq 0. \end{aligned} \tag{2.3}$$

More precisely, in limit (2.3), there are only two cases with respect to $u \to 0$, which are $u \downarrow 0$ and $u \uparrow 0$. We consider case by case below. By $x^* \in [-y^*, y^*]$ and $y^* \geq 0$, we have

$$\limsup_{u \downarrow 0} \frac{x^* u - y^*|u|}{2|u|} = \frac{x^* - y^*}{2} \leq 0. \tag{2.4}$$

And $$\limsup_{u\uparrow 0}\frac{x^*u\ - y^*|u|}{2|u|}=\frac{-x^*-y^*}{2}\leq 0. \qquad (2.5)$$

Substituting (2.4) and (2.5) into (2.3), we obtain again

$$\limsup_{u\to 0}\frac{\langle x^*,u-0\rangle_X\ -\ \langle y^*,f(u)-0\rangle_Y}{\|u-0\|_X+\|f(u)-0\|_Y}\leq 0.$$

For the given $y^*\geq 0$, this implies that

$$x^*\in\widehat{D}^*f(0,0)(y^*),\text{ for any } x^*\in[-y^*,y^*].$$

Hence, we proved that

$$[-y^*,y^*]\subseteq\widehat{D}^*f(0,0)(y^*),\text{ for any } y^*\geq 0. \qquad (2.6)$$

Let $x^*>\ y^*$. We estimate

$$\limsup_{u\to 0}\frac{\langle x^*,u-0\rangle_X\ -\ \langle y^*,f(u)-0\rangle_Y}{\|u-0\|_X+\|f(u)-0\|_Y}\geq\limsup_{u\downarrow 0}\frac{x^*u\ - y^*|u|}{2|u|}\ =\frac{x^*-y^*}{2}>0.$$

This implies that

$$x^*\notin\widehat{D}^*f(0,0)(y^*),\text{ for any } x^*>\ y^*. \qquad (2.7)$$

Next, let $x^*<-y^*$. We estimate

$$\limsup_{u\to 0}\frac{\langle x^*,u-0\rangle_X\ -\ \langle y^*,f(u)-0\rangle_Y}{\|u-0\|_X+\|f(u)-0\|_Y}\geq\limsup_{u\uparrow 0}\frac{x^*u\ - y^*|u|}{2|u|}=\frac{-x^*-y^*}{2}>0.$$

This implies that

$$x^*\notin\widehat{D}^*f(0,0)(y^*),\text{ for any } x^*<-y^*. \qquad (2.8)$$

(2.2) is proved by (2.6), (2.7) and (2.8). □

Proof of (c). Let $y^*<0$. There are two cases for consideration.

Case 1. Let $x^*\geq 0$. We calculate

$$\limsup_{u\to 0}\frac{\langle x^*,u-0\rangle_X\ -\ \langle y^*,f(u)-0\rangle_Y}{\|u-0\|_X+\|f(u)-0\|_Y}$$

$$\geq\limsup_{u\downarrow 0}\frac{x^*u\ - y^*|u|}{2|u|}$$

$$=\frac{x^*-y^*}{2}$$

$$>0.$$

Then, we obtain that, for $y^*<0$,

$$x^*\geq 0\Longrightarrow x^*\notin\widehat{D}^*f(0,0)(y^*).$$

Case 2. Let $x^* < 0$. We calculate

$$\limsup_{u\to 0} \frac{\langle x^*, u-0\rangle_X - \langle y^*, f(u)-0\rangle_Y}{\|u-0\|_X + \|f(u)-0\|_Y}$$

$$\geq \limsup_{u\uparrow 0} \frac{x^*u - y^*|u|}{2|u|}$$

$$= \frac{-x^*-y^*}{2}$$

$$> 0.$$

Then, we obtain that, for $y^* < 0$,

$$x^* < 0 \Longrightarrow x^* \notin \widehat{D}^*f(0,0)(y^*).$$

These prove (c). □

Under the sense of Mordukhovich differentiability defined by (1.3), part (b) and part (a) of Example 2.1 demonstrates (1.7). That is,

$f$ is Mordukhovich differentiable at $(0,0)$ and $f$ is not Fréchet differentiable at $0$.

Example 2.1 can be extended from $\mathbb{R}$ to $\mathbb{R}^n$. However, for simplicity, next, we extend it from $\mathbb{R}$ to $\mathbb{R}^2$.

**Example 2.2**. Let $X = X^* = Y = Y^* = \mathbb{R}^2$, in which the Hilbert norm in $\mathbb{R}^2$ is denoted by $\|\cdot\|_2$. Let $\langle\cdot,\cdot\rangle$ denote the inner product in $\mathbb{R}^2$. Define $g: \mathbb{R}^2 \to \mathbb{R}^2$, by

$$g(x_1, x_2) = (|x_1|, |x_2|), \text{ for any } (x_1, x_2) \in \mathbb{R}^2.$$

Then, by $g(0,0) = (0,0)$, the Fréchet differentiability of $g$ at point $(0,0)$ and the Mordukhovich differentiability of $g$ at point $((0,0), (0,0))$ have the following properties.

(a) $g$ is not Fréchet differentiable at $(0,0)$, that is,

$$\nabla g(0,0) \text{ does not exist.}$$

(b) For any $y^* = (y_1^*, y_2^*) \in \mathbb{R}^2$ with $y_i^* \geq 0$, for $i$ = 1, 2,

$$\widehat{D}^*g(0,0)(y^*) = \{x^* = (x_1^*, x_2^*) \in \mathbb{R}^2: |x_i^*| \leq y_i^*, \text{for } i = 1, 2\}. \tag{2.9}$$

(c) Let $y^* = (y_1^*, y_2^*) \in \mathbb{R}^2$. If $\min\{y_1^*, y_2^*\} < 0$, then

$$\widehat{D}^*g(0,0)(y^*) = \emptyset.$$

*Proof*. The proof of this example is similar to the proof of Example 2.1 that is an extension from $\mathbb{R}$ to $\mathbb{R}^2$.

Proof of (a). Assume that $\nabla g(0,0)$ exists, which is a continuous and linear operator from $\mathbb{R}^2$ to $\mathbb{R}^2$. So, $\nabla g(0,0)$ is represented by a real $2 \times 2$ matrix. Suppose

$$\nabla g(0,0) = \begin{pmatrix} a & b \\ c & d \end{pmatrix}, \text{ for some real numbers } a, b, c, d. \tag{2.10}$$

Under the assumption that $\nabla g(0,0)$ exists, by (1.3) and (2.10), we have

$$(0,0) = \lim_{(u_1,u_2)\to(0,0)} \frac{g\big((0,0)+(u_1,u_2)\big)-g(0,0)-\nabla g(0,0)(u_1,u_2)}{\|(u_1,u_2)\|_2}$$

$$= \lim_{(u_1,u_2)\to(0,0)} \frac{g(u_1,u_2)-\nabla g(0,0)(u_1,u_2)}{\|(u_1,u_2)\|_2}$$

$$= \lim_{(u_1,u_2)\to(0,0)} \frac{(|u_1|,|u_2|)-(u_1,u_2)\begin{pmatrix} a & b \\ c & d \end{pmatrix}}{\|(u_1,u_2)\|_2}$$

$$= \lim_{(u_1,u_2)\to(0,0)} \frac{(|u_1|,|u_2|)-(au_1+cu_2,\ bu_1+du_2)}{\|(u_1,u_2)\|_2}. \tag{2.11}$$

In particular, in (2.11), we have

$$(0,0) = \lim_{(u_1,0)\to(0,0),u_1\downarrow 0} \frac{(u_1,0)-(au_1,\ bu_1)}{u_1} = (1-a,\ -b);$$

$$(0,0) = \lim_{(u_1,0)\to(0,0),u_1\uparrow 0} \frac{(-u_1,0)-(au_1,\ bu_1)}{-u_1} = (1+a,\ b).$$

By adding the above two equations side by side, it gets $(0,0) = (2,0)$. This contradiction implies that $\nabla g(0,0)$ does not exist.

Proof of (b). Let $y^* = (y_1^*, y_2^*) \in \mathbb{R}^2$ with $y_i^* \geq 0$, for *i* = 1, 2. For $x^* = (x_1^*, x_2^*) \in \mathbb{R}^2$ with $|x_i^*| \leq y_i^*$, for *i* = 1, 2. We calculate

$$\limsup_{(u_1,u_2)\to(0,0)} \frac{\langle (x_1^*,x_2^*),\ (u_1,u_2)-(0,0)\rangle - \langle (y_1^*,y_2^*), g(u_1,u_2)-g(0,0)\rangle}{\|(u_1,u_2)-(0,0)\|_2 + \|g(u_1,u_2)-g(0,0)\|_2}$$

$$= \limsup_{(u_1,u_2)\to(0,0)} \frac{\langle (x_1^*,x_2^*),\ (u_1,u_2)\rangle - \langle (y_1^*,y_2^*),\ (|u_1|,|u_2|)\rangle}{\|(u_1,u_2)\|_2 + \|(|u_1|,|u_2|)\|_2}$$

$$= \limsup_{(u_1,u_2)\to(0,0)} \frac{x_1^*u_1+x_2^*u_2 - (y_1^*|u_1|+y_2^*|u_1|)}{2\|(u_1,u_2)\|_2}$$

$$= \limsup_{(u_1,u_2)\to(0,0)} \frac{x_1^*u_1-y_1^*|u_1| + x_2^*u_2 - y_2^*|u_1|}{2\|(u_1,u_2)\|_2}$$

$$\leq 0.$$

This implies that

$$x^* \in \widehat{D}^* g(0,0)(y^*), \text{ for any } x^* = (x_1^*, x_2^*) \in \mathbb{R}^2 \text{ with } |x_i^*| \leq y_i^*, \text{ for } i = 1, 2. \tag{2.12}$$

Let $x^* = (x_1^*, x_2^*) \in \mathbb{R}^2$. Suppose that the inequalities $|x_i^*| \leq y_i^*$, for *i* = 1, 2 do not hold. Without loose the generality, assume that $x_1^* > y_1^*$. We calculate

$$\limsup_{(u_1,u_2)\to(0,0)} \frac{\langle (x_1^*,x_2^*),\ (u_1,u_2)-(0,0)\rangle - \langle (y_1^*,y_2^*), g(u_1,u_2)-g(0,0)\rangle}{\|(u_1,u_2)-(0,0)\|_2 + \|g(u_1,u_2)-g(0,0)\|_2}$$

$$\geq \limsup_{(u_1,0)\to(0,0),u_1\downarrow 0} \frac{\langle (x_1^*,x_2^*),\ (u_1,0)\rangle - \langle (y_1^*,y_2^*),\ (|u_1|,|0|)\rangle}{\|(u_1,0)\|_2 + \|(|u_1|,|0|)\|_2}$$

$$= \limsup_{(u_1,0)\to(0,0),u_1\downarrow 0} \frac{x_1^* u_1 - y_1^* u_1}{2u_1}$$

$$= \frac{x_1^* - y_1^*}{2}$$

$$> 0.$$

We obtain that

$$x_1^* > y_1^* \Longrightarrow x^* \notin \widehat{D}^* g(0,0)(y^*). \tag{2.13}$$

Let $x^* = (x_1^*, x_2^*) \in \mathbb{R}^2$. Without loose the generality, suppose $x_1^* < -y_1^*$. We calculate

$$\limsup_{(u_1,u_2)\to(0,0)} \frac{\langle (x_1^*,x_2^*),\ (u_1,u_2)-(0,0)\rangle - \langle (y_1^*,y_2^*), g(u_1,u_2)-g(0,0)\rangle}{\|(u_1,u_2)-(0,0)\|_2 + \|g(u_1,u_2)-g(0,0)\|_2}$$

$$\geq \limsup_{(u_1,0)\to(0,0),u_1\uparrow 0} \frac{\langle (x_1^*,x_2^*),\ (u_1,0)\rangle - \langle (y_1^*,y_2^*),\ (|u_1|,|0|)\rangle}{\|(u_1,0)\|_2 + \|(|u_1|,|0|)\|_2}$$

$$= \limsup_{(u_1,0)\to(0,0),u_1\uparrow 0} \frac{x_1^* u_1 - y_1^* |u_1|}{2|u_1|}$$

$$= \frac{-x_1^* - y_1^*}{2}$$

$$> 0.$$

We obtain that

$$x_1^* < -y_1^* \Longrightarrow x^* \notin \widehat{D}^* g(0,0)(y^*). \tag{2.14}$$

Then, (2.9) in part (b) is proved by (2.12), (2.13) and (2.14).

Proof of (c). Let $y^* = (y_1^*, y_2^*) \in \mathbb{R}^2$ with $\min\{y_1^*, y_2^*\} < 0$. Without loose the generality, suppose $y_1^* < 0$. Let arbitrarily given $x^* = (x_1^*, x_2^*) \in \mathbb{R}^2$. We consider two cases:

Case 1. Suppose $x_1^* \geq 0$. By $y_1^* < 0$, we calculate

$$\limsup_{(u_1,u_2)\to(0,0)} \frac{\langle (x_1^*,x_2^*),\ (u_1,u_2)-(0,0)\rangle - \langle (y_1^*,y_2^*), g(u_1,u_2)-g(0,0)\rangle}{\|(u_1,u_2)-(0,0)\|_2 + \|g(u_1,u_2)-g(0,0)\|_2}$$

$$\geq \limsup_{(u_1,0)\to(0,0),u_1\downarrow 0} \frac{\langle (x_1^*,x_2^*),\ (u_1,0)\rangle - \langle (y_1^*,y_2^*),\ (|u_1|,|0|)\rangle}{\|(u_1,0)\|_2 + \|(|u_1|,|0|)\|_2}$$

$$= \limsup_{(u_1,0)\to(0,0),u_1\downarrow 0} \frac{x_1^* u_1 - y_1^* |u_1|}{2|u_1|}$$

$$= \frac{x_1^* - y_1^*}{2}$$

$$> 0.$$

Then, for given $x^* = (x_1^*, x_2^*) \in \mathbb{R}^2$, we obtain that

$$x_1^* \geq 0 \implies x^* \notin \widehat{D}^* g(0,0)(y^*). \tag{2.15}$$

Case 2. Suppose $x_1^* < 0$. By $y_1^* < 0$, we calculate

$$\limsup_{(u_1,u_2)\to(0,0)} \frac{\langle (x_1^*,x_2^*),\ (u_1,u_2)-(0,0)\rangle - \langle (y_1^*,y_2^*), g(u_1,u_2)-g(0,0)\rangle}{\|(u_1,u_2)-(0,0)\|_2 + \|g(u_1,u_2)-g(0,0)\|_2}$$

$$\geq \limsup_{(u_1,0)\to(0,0), u_1\uparrow 0} \frac{\langle (x_1^*,x_2^*),\ (u_1,0)\rangle - \langle (y_1^*,y_2^*),\ (|u_1|,|0|)\rangle}{\|(u_1,0)\|_2 + \|(|u_1|,|0|)\|_2}$$

$$= \limsup_{(u_1,0)\to(0,0), u_1\uparrow 0} \frac{x_1^* u_1 - y_1^*|u_1|}{2|u_1|}$$

$$= \frac{-x_1^* - y_1^*}{2}$$

$$> 0.$$

Then, for given $x^* = (x_1^*, x_2^*) \in \mathbb{R}^2$, we obtain that

$$x_1^* < 0 \implies x^* \notin \widehat{D}^* g(0,0)(y^*). \tag{2.16}$$

By (2.15) and (2.16), part (c) is proved. That is, for any $y^* = (y_1^*, y_2^*) \in \mathbb{R}^2$

$$\min\{y_1^*, y_2^*\} < 0 \implies \widehat{D}^* g(0,0)(y^*) = \emptyset. \qquad \square$$

Once again, under the sense of Mordukhovich differentiability defined by (1.3), part (b) and part (a) of Example 2.2 demonstrate (1.7). That is,

$g$ is Mordukhovich differentiable at $((0,0),(0,0))$ and $g$ is not Fréchet differentiable at $(0,0)$.

In the following example, we construct a single-valued mapping $h: \mathbb{R} \to \mathbb{R}$ that satisfies that $h$ is nowhere Fréchet differentiable and it is Mordukhovich differentiable at every point.

**Example 2.3**. Define a single-valued mapping $h: \mathbb{R} \to \mathbb{R}$ by

$$h(x) = \begin{cases} 1, & \text{if } x \text{ is rational,} \\ -1, & \text{if } x \text{ is irrational.} \end{cases}$$

Let $(\bar{x}, \bar{y}) \in \mathcal{G}(h)$. The Fréchet differentiability of $h$ at point $\bar{x}$ and the Mordukhovich differentiability of $h$ at point $(\bar{x}, \bar{y})$ have the following properties.

(i) $h$ is nowhere Fréchet differentiable, that is,

$$\nabla h(\bar{x}) \text{ does not exist, for any } \bar{x} \in \mathbb{R};$$

(ii) The Mordukhovich derivative $\widehat{D}^* h(\bar{x}, \bar{y})$ of $h$ satisfies that

$$\widehat{D}^* h(\bar{x}, \bar{y})(y^*) \cap (\mathbb{R}\backslash\{0\}) = \emptyset, \text{ for any } y^* \in \mathbb{R}. \tag{2.17}$$

However, $h$ is Mordukhovich differentiable at every $(\bar{x}, \bar{y}) \in \mathcal{G}(h)$. More precisely, we have that

(a) if $\bar{x}$ is rational then, $(\bar{x}, \bar{y}) = (\bar{x}, 1) \in \mathcal{G}(h)$ and

$$\widehat{D}^*h(\bar{x},1)(y^*) = \{0\}, \text{ for any } y^* \leq 0;$$

(b) if $\bar{x}$ is irrational, then, $(\bar{x},\bar{y}) = (\bar{x},-1) \in \mathcal{G}(h)$ and

$$\widehat{D}^*h(\bar{x},-1)(y^*) = \{0\}, \text{ for any } y^* \leq 0.$$

*Proof*. The proof of (i) is straight forward and it is omitted here.

Proof of (ii). There are two cases for $\bar{x}$ to be rational or irrational.

Case 1. $(\bar{x},\bar{y}) \in \mathcal{G}(h)$ with $\bar{x}$ being rational. By definition, $\bar{y} = h(\bar{x}) = 1$. Let $y^* \in \mathbb{R}$ be arbitrarily given. For any $x^* \in \mathbb{R}\backslash\{0\}$, we calculate

$$\limsup_{\substack{(u,v)\to(\bar{x},1)\\ v=h(u)}} \frac{\langle x^*,u-\bar{x}\rangle - \langle y^*,v-1\rangle}{|u-\bar{x}| + |v-1|}$$

$$\geq \limsup_{\substack{(u,\bar{y})\to(\bar{x},\bar{y})\\ u \text{ is rational and } 1=h(u)}} \frac{\langle x^*,u-\bar{x}\rangle - \langle y^*,1-1\rangle}{|u-\bar{x}| + |1-1|}$$

$$= \begin{cases} x^*, & \text{for } u \downarrow \bar{x} \text{ and } u \text{ is rational,} \\ -x^*, & \text{for } u \uparrow \bar{x} \text{ and } u \text{ is rational.} \end{cases}$$

This implies that

$$\limsup_{\substack{(u,v)\to(\bar{x},1)\\ v=h(u)}} \frac{\langle x^*,u-\bar{x}\rangle - \langle y^*,v-1\rangle}{|u-\bar{x}| + |v-1|} \geq |x^*| > 0.$$

We obtain that, for any $y^* \in \mathbb{R}$,

$$x^* \neq 0 \quad \Longrightarrow \quad x^* \notin \widehat{D}^*F(\bar{x},\bar{y})(y^*), \text{ for any } y^* \in \mathbb{R}. \tag{2.18}$$

Case 2. $(\bar{x},\bar{y}) \in \mathcal{G}(h)$ with $\bar{x}$ being irrational. By definition, $\bar{y} = h(\bar{x}) = -1$. In this case, (3.13) can be similarly proved. This proves part (2.17).

Proof of (a). In particular, if $\bar{x}$ is rational and $\bar{y} = h(\bar{x}) = 1$, then, for any $y^* \leq 0$,

$$\limsup_{\substack{(u,v)\to(\bar{x},1)\\ v=h(u)}} \frac{\langle 0,u-\bar{x}\rangle - \langle y^*,v-1\rangle}{|u-\bar{x}| + |v-1|} = \limsup_{\substack{(u,v)\to(\bar{x},1)\\ v=h(u)}} \frac{-\langle y^*,v-1\rangle}{|u-\bar{x}| + |v-1|} \leq 0.$$

By (2.18), this proves part (a).

Proof of (b). In particular, if $\bar{x}$ is irrational and $\bar{y} = h(\bar{x}) = -1$, then, for any $y^* \geq 0$,

$$\limsup_{\substack{(u,v)\to(\bar{x},-1)\\ v=h(u)}} \frac{\langle 0,u-\bar{x}\rangle - \langle y^*,v+1\rangle}{|u-\bar{x}| + |v+1|} = \limsup_{\substack{(u,v)\to(\bar{x},1)\\ v=h(u)}} \frac{-\langle y^*,v+1\rangle}{|u-\bar{x}| + |v+1|} \leq 0.$$

By (2.18), this proves part (b). □

The results of Example 2.3 are stronger than examples 2.1 and 2.2. Under the sense of Mordukhovich differentiability defined by (1.3), part (b) and part (a) of Example 2.3 demonstrate (1.7). That is, for any $\bar{x} \in \mathbb{R}$ with $(\bar{x}, \bar{y}) \in \mathcal{G}(h)$,

$$h \text{ is Mordukhovich differentiable at } (\bar{x}, \bar{y}) \quad \text{and} \quad h \text{ is not Fréchet differentiable at } \bar{x}.$$

## 3. Some Properties of Mordukhovich Derivatives of Set-valued Mappings

Recall that in this paper, $(X, \|\cdot\|_X)$ and $(Y, \|\cdot\|_Y)$ are real Banach spaces with topological dual spaces $(X^*, \|\cdot\|_{X^*})$ and $(Y^*, \|\cdot\|_{Y^*})$, respectively. Let $\langle\cdot, \cdot\rangle_X$ denote the real canonical pairing between $X^*$ and $X$ and $\langle\cdot, \cdot\rangle_Y$ the real canonical pairing between $Y^*$ and $Y$. Let $\theta_X$ and $\theta_Y$ denote the origins in $X$ and $Y$, respectively. Let $\theta_X^*$ and $\theta_Y^*$ denote the null elements in $X^*$ and $Y^*$, respectively. Let $A$ be a nonempty convex open subset of $X$.

**Theorem 3.1**. *Let $F: A \rightrightarrows Y$ be a set-valued mapping. Let $C$ be a nonempty convex open subset of $A$. Let $\bar{y} \in Y$ satisfying*

$$\bar{y} \in F(x), \text{ for each } x \in C.$$

*Then, for any $\bar{x} \in C$, the Mordukhovich derivative $\widehat{D}^*F(\bar{x}, \bar{y}): Y^* \rightrightarrows X^*$ has the following properties*:

(a) $\widehat{D}^*F(\bar{x}, \bar{y})(\theta_Y^*) = \{\theta_X^*\}$;
(b) *For any $y^* \in Y^* \setminus \{\theta_Y^*\}$,*

$$x^* \notin \widehat{D}^*F(\bar{x}, \bar{y})(y^*), \text{ for any } x^* \in X^* \setminus \{\theta_X^*\}. \tag{3.1}$$

*Proof*. Proof of (a). It is clear that

$$\theta_X^* \in \widehat{D}^*F(\bar{x}, \bar{y})(\theta_Y^*). \tag{3.2}$$

Next, for any $x^* \in X^* \setminus \{\theta_X^*\}$, there is $w \in X \setminus \{\theta_X\}$ such that $\langle x^*, w\rangle_X > 0$. For any real number $t$, let $u(t) = \bar{x} + tw$. Then, there is $\delta > 0$, such that

$$u(t) = \bar{x} + tw \in C, \text{ for any } |t| < \delta.$$

We calculate

$$\limsup_{\substack{(u,v)\to(\bar{x},\bar{y}) \\ u\in A \text{ and } v\in F(u)}} \frac{\langle x^*, u-\bar{x}\rangle_X - \langle \theta_Y^*, v-\bar{y}\rangle_Y}{\|u-\bar{x}\|_X + \|v-\bar{x}\|_Y}$$

$$\geq \limsup_{\substack{t\to 0 \text{ with } |t|<\delta \\ \bar{y}\in F(u(t))}} \frac{\langle x^*, u(t)-\bar{x}\rangle_X}{\|u(t)-\bar{x}\|_X + \|\bar{y}-\bar{y}\|_Y}$$

$$= \limsup_{\substack{t\to 0 \text{ with } |t|<\delta \\ \bar{y}\in F(u(t))}} \frac{\langle x^*, \bar{x} + tw - \bar{x}\rangle_X}{\|\bar{x} + tw - \bar{x}\|_X}$$

$$\geq \limsup_{\substack{t\downarrow 0 \text{ with } 0 < t<\delta \\ \bar{y}\in F(u(t))}} \frac{\langle x^*, tw\rangle_X}{\|tw\|_X}$$

$$= \limsup_{\substack{t\downarrow 0 \text{ with } 0<t<\delta \\ \bar{y}\in F(u(t))}} \frac{t\,\langle x^*,\ w\rangle_X}{t\|w\|_X}$$

$$= \frac{\langle x^*,\ w\rangle_X}{\|w\|_X}$$

$$> 0.$$

This implies that

$$x^* \notin \widehat{D}^*F(\bar{x},\bar{y})(\theta_Y^*), \text{ for any } x^* \in X^*\backslash\{\theta_X^*\}. \tag{3.3}$$

By (3.2) and (3.3), we obtain that $\widehat{D}^*F(\bar{x},\bar{y})(\theta_Y^*) = \{\theta_X^*\}$, which proves (a).

Proof of (b). Let $y^* \in Y^*\backslash\{\theta_Y^*\}$. For any $x^* \in X^*\backslash\{\theta_X^*\}$, similarly to the proof of (a), there is $w \in X\backslash\{\theta_X\}$ such that $\langle x^*, w\rangle_X > 0$. For any real number $t$, let $u(t) = \bar{x} + tw$. Then, there is $\delta > 0$, such that

$$u(t) = \bar{x} + tw \in C, \text{ for any } |t| < \delta.$$

We calculate

$$\limsup_{\substack{(u,v)\to(\bar{x},\bar{y}) \\ u\in A \text{ and } v\in F(u)}} \frac{\langle x^*,u-\bar{x}\rangle_X - \langle y^*,v-\bar{y}\rangle_Y}{\|u-\bar{x}\|_X + \|v-\bar{x}\|_Y}$$

$$\geq \limsup_{\substack{t\downarrow 0 \text{ with } 0<t<\delta \\ \bar{y}\in F(u(t))}} \frac{\langle x^*,u(t)-\bar{x}\rangle_X - \langle y^*,\bar{y}-\bar{y}\rangle_Y}{\|u(t)-\bar{x}\|_X + \|\bar{y}-\bar{y}\|_Y}$$

$$= \limsup_{\substack{t\downarrow 0 \text{ with } 0<t<\delta \\ \bar{y}\in F(u(t))}} \frac{\langle x^*,\ \bar{x}+tw-\bar{x}\rangle_X}{\|\bar{x}+tw-\bar{x}\|_X}$$

$$= \frac{\langle x^*,\ w\rangle_X}{\|w\|_X}$$

$$> 0.$$

This implies that

$$x^* \notin \widehat{D}^*F(\bar{x},\bar{y})(y^*), \text{ for any } x^* \in X^*\backslash\{\theta_X^*\}.$$

This proves (b). □

**Theorem 3.2**. *Let $F: A \rightrightarrows Y$ be a set-valued mapping. Suppose $\mathcal{G}(F)^o \neq \emptyset$. Then, for any $(\bar{x},\bar{y}) \in \mathcal{G}(F)^o$, the Mordukhovich derivative $\widehat{D}^*F(\bar{x},\bar{y}): Y^* \rightrightarrows X^*$ has the following properties*:

(a) $\widehat{D}^*F(\bar{x},\bar{y})(\theta_Y^*) = \{\theta_X^*\}$;
(b) *For any $y^* \in Y^*\backslash\{\theta_Y^*\}$,*

$$\widehat{D}^*F(\bar{x},\bar{y})(y^*) = \emptyset. \tag{3.4}$$

*Here, $\mathcal{G}(F)^o$ is the topological interior of $\mathcal{G}(F)$ with respect to the product topology of the topologies of X and Y.*

*Proof*. Proof of (a). Let $(\bar{x},\bar{y}) \in \mathcal{G}(F)^o$. There are open neighborhoods $C$ of $\bar{x}$ in $X$ and $D$ of $\bar{y}$ in $Y$ such that $C\times D \subseteq \mathcal{G}(F)$, which implies that

$$\{(x,\bar{y}) \in X \times Y: x \in C\} \subseteq C \times D \subseteq \mathcal{G}(F).$$

By $(\bar{x},\bar{y}) \in \{(x,\bar{y}) \in X \times Y: x \in C\}$, part (a) of this theorem follows from **Theorem 3.1** immediately. Then, we prove (b) of this theorem. At first, we prove that, for any $y^* \in Y^*\backslash\{\theta_Y^*\}$,

$$\theta_X^* \notin \widehat{D}^*F(\bar{x},\bar{y})(y^*). \tag{3.5}$$

For arbitrarily given $y^* \in Y^*\backslash\{\theta_Y^*\}$, there is $y \in Y$ such that $\langle y^*, y\rangle_Y < 0$. Let $z \in X\backslash\{\theta_X\}$. For any real number $s$, let $u(s) = \bar{x} + sz$. By the conditions of $C$ and $D$, there is $\gamma > 0$ such that

$$\bar{y} + sy \in D \text{ and } u(s) = \bar{x} + sz \in C, \text{ for } |s| < \gamma.$$

This implies that

$$\{(x,\bar{y} + sy) \in X \times Y: x \in C, |s| < \gamma\} \subseteq C \times D \subseteq \mathcal{G}(F).$$

In particular, by the definition of $\mathcal{G}(F)$, we have

$$\bar{y} + sy \in F(u), \text{ for any } u \in C \text{ and for any } |s| < \gamma.$$

Then, we calculate

$$\limsup_{\substack{(u,v)\to(\bar{x},\bar{y})\\ u\in A \text{ and } v\in F(u)}} \frac{\langle\theta_X^*,u-\bar{x}\rangle_X - \langle y^*,v-\bar{y}\rangle_Y}{\|u-x\|_X + \|v-\bar{x}\|_Y}$$

$$\geq \limsup_{\substack{s\to 0 \text{ with } |s|<\gamma\\ \bar{y}+sy \in F(u(s))}} \frac{-\langle y^*,\bar{y}+sy-\bar{y}\rangle_Y}{\|\bar{x} + sz-\bar{x}\|_X + \|\bar{y}+sy-\bar{y}\|_Y}$$

$$\geq \limsup_{\substack{s\downarrow 0 \text{ with } |s|<\gamma\\ \bar{y}+sy \in F(u(s))}} \frac{-s\langle y^*,y\rangle_Y}{s\|z\|_X+s\|y\|_Y}$$

$$= \frac{-\langle y^*,y\rangle_Y}{\|z\|_X+\|y\|_Y}$$

$$> 0.$$

This proves (3.5). Next, we show that for any given $y^* \in Y^*\backslash\{\theta_Y^*\}$,

$$x^* \in X^*\backslash\{\theta_X^*\} \implies x^* \notin \widehat{D}^*F(\bar{x},\bar{y})(y^*). \tag{3.6}$$

For any given $x^* \in X^*\backslash\{\theta_X^*\}$, there is $w \in X\backslash\{\theta_X\}$ such that $\langle x^*, w\rangle_X > 0$. For any real number $t$, let $u(t) = \bar{x} + tw$. Then, there is $\delta > 0$, such that

$$u(t) = \bar{x} + tw \in C, \text{ for any } |t| < \delta.$$

By $C\times D \subseteq \mathcal{G}(F)$, by the definition of $\mathcal{G}(F)$ and $\bar{y} \in D$, this implies that

$$\bar{y} \in F(u(t)), \text{ for any } |t| < \delta.$$

Similarly to the proof of Theorem 3.1, we calculate

$$\limsup_{\substack{(u,v)\to(\bar{x},\bar{y})\\ u\in A \text{ and } v\in F(u)}} \frac{\langle x^*, u-\bar{x}\rangle_X - \langle y^*, v-\bar{y}\rangle_Y}{\|u-\bar{x}\|_X + \|v-\bar{x}\|_Y}$$

$$\geq \limsup_{\substack{t\to 0 \text{ with } |t|<\delta\\ \bar{y}\in F(u(t))}} \frac{\langle x^*, u(t)-\bar{x}\rangle_X - \langle y^*, \bar{y}-\bar{y}\rangle_Y}{\|u(t)-\bar{x}\|_X + \|\bar{y}-\bar{y}\|_Y}$$

$$= \limsup_{\substack{t\to 0 \text{ with } |t|<\delta\\ \bar{y}\in F(u(t))}} \frac{\langle x^*, \ \bar{x}+tw-\bar{x}\rangle_X}{\|\bar{x}+tw-\bar{x}\|_X}$$

$$\geq \limsup_{\substack{t\downarrow 0 \text{ with } |t|<\delta\\ \bar{y}\in F(u(t))}} \frac{\langle x^*, \ tw\rangle_X}{\|tw\|_X}$$

$$= \limsup_{\substack{t\downarrow 0 \text{ with } |t|<\delta\\ \bar{y}\in F(u(t))}} \frac{t\,\langle x^*, \ w\rangle_X}{t\|w\|_X}$$

$$= \frac{\langle x^*, \ w\rangle_X}{\|w\|_X}$$

$$> 0.$$

This implies (3.6). Then, (3.4) is proved by (3.5) and (3.6). □

Next, we give some examples of set-valued mappings to demonstrate Theorem 3.2.

**Example 3.3**. Let $X = \mathbb{R}$, $Y = \mathbb{R}^2$ with null element $\theta = (0,0)$. Recall that $\|\cdot\|_2$ denote the Hilbert norm on $\mathbb{R}^2$. Define $F: \mathbb{R} \rightrightarrows \mathbb{R}^2$, by

$$F(x) = \{(y_1, y_2) \in \mathbb{R}^2 : y_1^2 + y_1^2 \leq x^2\}, \text{ for any } x \in \mathbb{R}. \tag{3.7}$$

$F$ is a set-valued mapping with closed and convex values. The graph of $F$ consists of two apposite regular cones in $\mathbb{R}^3$, which share the same vertex (0, 0, 0), and the $x$-axis is the center axis of these regular cones, which satisfies the following equation.

$$\mathcal{G}(F) = \{(x, y_1, y_2) \in \mathbb{R}^3 : y_1^2 + y_1^2 \leq x^2, x \in \mathbb{R}\}.$$

It is clear that the interior of $\mathcal{G}(F)$ satisfies the following equation.

$$\mathcal{G}(F)^o = \{(x, y_1, y_2) \in \mathbb{R}^3 : y_1^2 + y_1^2 < x^2, x \in \mathbb{R}\backslash\{0\}\}. \tag{3.8}$$

The Mordukhovich differentiability of $F$ has the following properties.

(i) For any $(\bar{x}, \bar{y}) = (\bar{x}, (\bar{y}_1, \bar{y}_2)) \in \mathcal{G}(F)^o$ with $\bar{x} \in \mathbb{R}\backslash\{0\}$, $\bar{y} = (\bar{y}_1, \bar{y}_2) \in \mathbb{R}^2$ satisfying $\bar{y}_1^2 + \bar{y}_1^2 < \bar{x}^2$, we have

(a) $\widehat{D}^* F(\bar{x}, \bar{y})(\theta) = \{0\}$;

(b) For any $y^* \in \mathbb{R}^2\backslash\{\theta\}$,

$$\widehat{D}^*F(\bar{x},\bar{y})(y^*) = \emptyset. \tag{3.9}$$

(ii) Let $(\bar{x},\bar{y}_1,\bar{y}_2) \in \partial\mathcal{G}(F)$ with $\bar{x} \neq 0$, $\bar{y} = (\bar{y}_1,\bar{y}_2) \in \mathbb{R}^2$ satisfying $\bar{y}_1^2 + \bar{y}_1^2 = \bar{x}^2$. Then, for any $y^* = (y_1^*, y_2^*) \in \mathbb{R}^2$,

(a) $$y_1^*\bar{y}_1 + y_2^*\bar{y}_2 \neq 0 \quad \Longrightarrow \quad 0 \notin \widehat{D}^*F(\bar{x},\bar{y})(y^*);$$

(b) $$\big((-\infty, -\|y^*\|_2) \cup (\|y^*\|_2, \infty)\big) \cap \widehat{D}^*F(\bar{x},\bar{y})(y^*) = \emptyset.$$

(iii) For any $y^* = (y_1^*, y_2^*) \in \mathbb{R}^2$,

$$x^* \neq 0 \quad \Longrightarrow \quad x^* \notin \widehat{D}^*F(0,\theta)(y^*).$$

Here, $\partial\mathcal{G}(F)$ denotes the boundary of $\mathcal{G}(F)$.

*Proof.* Proof of (i). Part (i) and (3.9) are proved by (3.7), (3.8) and Theorem 3.2 immediately.

Proof of (a) in (ii). In this part, we consider the Mordukhovich differentiability of $F$ on the boundary of $\mathcal{G}(F)$, which is denoted by $\partial\mathcal{G}(F)$.

Let $(\bar{x},\bar{y}) \in \partial\mathcal{G}(F)$ with $\bar{x} \neq 0$ and $\bar{y} = (\bar{y}_1,\bar{y}_2) \in \mathbb{R}^2$ satisfying $\bar{y}_1^2 + \bar{y}_1^2 = \bar{x}^2$. Let $y^* = (y_1^*, y_2^*) \in \mathbb{R}^2$ with $y_1^*\bar{y}_1 + y_2^*\bar{y}_2 \neq 0$ (which implies $(y_1^*, y_2^*) \neq \theta$). Let $x^* = 0$. By $\bar{x} \neq 0$ and $\bar{y}_1^2 + \bar{y}_1^2 = \bar{x}^2$, we calculate

$$\limsup_{\substack{(u,(v_1,v_2))\to(\bar{x},(\bar{y}_1,\bar{y}_2))\\ u\in\mathbb{R} \text{ and } (v_1,v_2)\in F(u)}} \frac{\langle x^*, u-\bar{x}\rangle - \langle (y_1^*,y_2^*),\ (v_1,v_2)-(\bar{y}_1,\bar{y}_2)\rangle}{|u-\bar{x}| + \|(v_1,v_2)-(\bar{y}_1,\bar{y}_2)\|_2}$$

$$= \limsup_{\substack{(u,(v_1,v_2))\to(\bar{x},(\bar{y}_1,\bar{y}_2))\\ u\in\mathbb{R} \text{ and } (v_1,v_2)\in F(u)}} \frac{x^*(u-\bar{x}) - (y_1^*(v_1-\bar{y}_1) + y_2^*(v_2-\bar{y}_2))}{|u-\bar{x}| + \|(v_1,v_2)-(\bar{y}_1,\bar{y}_2)\|_2}$$

$$= \limsup_{\substack{(u,(v_1,v_2))\to(\bar{x},(\bar{y}_1,\bar{y}_2))\\ u\in\mathbb{R} \text{ and } (v_1,v_2)\in F(u)}} \frac{-(y_1^*(v_1-\bar{y}_1) + y_2^*(v_2-\bar{y}_2))}{|u-\bar{x}| + \|(v_1,v_2)-(\bar{y}_1,\bar{y}_2)\|_2}$$

$$\geq \limsup_{\substack{t(\bar{x},(\bar{y}_1,\bar{y}_2))\to(\bar{x},(\bar{y}_1,\bar{y}_2))\\ t\to 1}} \frac{-(t-1)(y_1^*\bar{y}_1 + y_2^*\bar{y}_2)}{|(t-1)\bar{x}| + \|(t-1)(\bar{y}_1,\bar{y}_2)\|_2}$$

$$\geq \limsup_{\substack{t(\bar{x},(\bar{y}_1,\bar{y}_2))\to(\bar{x},(\bar{y}_1,\bar{y}_2))\\ t\to 1}} \frac{-(t-1)(y_1^*\bar{y}_1 + y_2^*\bar{y}_2)}{|t-1||\bar{x}| + |t-1||\bar{x}|} \quad (\bar{y}_1^2 + \bar{y}_1^2 = \bar{x}^2)$$

$$\geq \limsup_{\substack{t(\bar{x},(\bar{y}_1,\bar{y}_2))\to(\bar{x},(\bar{y}_1,\bar{y}_2))\\ t\to 1}} \frac{-(t-1)(y_1^*\bar{y}_1 + y_2^*\bar{y}_2)}{2|t-1||\bar{x}|}$$

$$= \begin{cases} \dfrac{-(y_1^*\bar{y}_1 + y_2^*\bar{y}_2)}{2|\bar{x}|}, & \text{for } t \downarrow 1, \\[2ex] \dfrac{y_1^*\bar{y}_1 + y_2^*\bar{y}_2}{2|\bar{x}|}, & \text{for } t \uparrow 1. \end{cases}$$

This implies that

$$\limsup_{\substack{(u,(v_1,v_2))\to(\bar{x},(\bar{y}_1,\bar{y}_2))\\ u\in\mathbb{R} \text{ and } (v_1,v_2)\in F(u)}} \frac{-(y_1^*(v_1-\bar{y}_1)+y_2^*(v_2-\bar{y}_2))}{|u-\bar{x}|+\|(v_1,v_2)-(\bar{y}_1,\bar{y}_2)\|_2} \geq \left|\frac{y_1^*\bar{y}_1+y_2^*\bar{y}_2}{2\bar{x}}\right| > 0.$$

Hence, we obtain that, for $(\bar{x},\bar{y}_1,\bar{y}_2)\in\partial\mathcal{G}(F)$ with $\bar{x}\neq 0$, $\bar{y}=(\bar{y}_1,\bar{y}_2)\in\mathbb{R}^2$ satisfying $\bar{y}_1^2+\bar{y}_1^2=\bar{x}^2$,

$$0\notin \widehat{D}^*F(\bar{x},(\bar{y}_1,\bar{y}_2))(y^*), \text{ for any } y^*=(y_1^*,y_2^*)\in\mathbb{R}^2 \text{ with } y_1^*\bar{y}_1+y_2^*\bar{y}_2\neq 0.$$

This proves (a) in (ii) of this example.

Proof of (b) in (ii). At first, we suppose $\bar{x}>0$. Then, the proof of (b) in (ii) is divided to two cases.

Case 1. Let $x^*\in(\|y^*\|_2,\infty)$. Similarly to the proof of (a) in (ii), by $\bar{x}>0$ and $\bar{y}_1^2+\bar{y}_1^2=\bar{x}^2$, we calculate

$$\begin{aligned}
&\limsup_{\substack{(u,(v_1,v_2))\to(\bar{x},(\bar{y}_1,\bar{y}_2))\\ u\in\mathbb{R} \text{ and } (v_1,v_2)\in F(u)}} \frac{\langle x^*,u-\bar{x}\rangle-\langle(y_1^*,y_2^*),\ (v_1,v_2)-(\bar{y}_1,\bar{y}_2)\rangle}{|u-\bar{x}|+\|(v_1,v_2)-(\bar{y}_1,\bar{y}_2)\|_2}\\
&=\limsup_{\substack{(u,(v_1,v_2))\to(\bar{x},(\bar{y}_1,\bar{y}_2))\\ u\in\mathbb{R} \text{ and } (v_1,v_2)\in F(u)}} \frac{x^*(u-\bar{x})-(y_1^*(v_1-\bar{y}_1)+y_2^*(v_2-\bar{y}_2))}{|u-\bar{x}|+\|(v_1,v_2)-(\bar{y}_1,\bar{y}_2)\|_2}\\
&\geq\limsup_{\substack{t(\bar{x},(\bar{y}_1,\bar{y}_2))\to(\bar{x},(\bar{y}_1,\bar{y}_2))\\ t\to 1}} \frac{(t-1)x^*\bar{x}-(t-1)(y_1^*\bar{y}_1+y_2^*\bar{y}_2)}{2|(t-1)\bar{x}|}\\
&\geq\limsup_{\substack{t(\bar{x},(\bar{y}_1,\bar{y}_2))\to(\bar{x},(\bar{y}_1,\bar{y}_2))\\ t\downarrow 1}} \frac{x^*\bar{x}-(y_1^*\bar{y}_1+y_2^*\bar{y}_2)}{2\bar{x}}\\
&\geq\limsup_{\substack{t(\bar{x},(\bar{y}_1,\bar{y}_2))\to(\bar{x},(\bar{y}_1,\bar{y}_2))\\ t\downarrow 1}} \frac{x^*\bar{x}-\|(\bar{y}_1,\bar{y}_2)\|_2\|(y_1^*,y_2^*)\|_2}{2\bar{x}}\\
&=\frac{x^*-\|(y_1^*,y_2^*)\|_2}{2} \quad (\text{by } \bar{x}>0 \text{ and } \bar{y}_1^2+\bar{y}_1^2=\bar{x}^2)\\
&\quad>0.
\end{aligned}$$

This implies that, for any $x^*\in\mathbb{R}$,

$$x^*\in(\|y^*\|_2,\infty) \quad\Longrightarrow\quad x^*\notin\widehat{D}^*F(\bar{x},(\bar{y}_1,\bar{y}_2))(y^*). \tag{3.10}$$

Case 2. Let $x^*\in(-\infty,-\|y^*\|_2)$. By $\bar{x}>0$ and $\bar{y}_1^2+\bar{y}_1^2=\bar{x}^2$, we calculate

$$\begin{aligned}
&\limsup_{\substack{(u,(v_1,v_2))\to(\bar{x},(\bar{y}_1,\bar{y}_2))\\ u\in\mathbb{R} \text{ and } (v_1,v_2)\in F(u)}} \frac{\langle x^*,u-\bar{x}\rangle-\langle(y_1^*,y_2^*),\ (v_1,v_2)-(\bar{y}_1,\bar{y}_2)\rangle}{|u-\bar{x}|+\|(v_1,v_2)-(\bar{y}_1,\bar{y}_2)\|_2}\\
&\geq\limsup_{\substack{t(\bar{x},(\bar{y}_1,\bar{y}_2))\to(\bar{x},(\bar{y}_1,\bar{y}_2))\\ t\to 1}} \frac{(t-1)x^*\bar{x}-(t-1)(y_1^*\bar{y}_1+y_2^*\bar{y}_2)}{2|(t-1)\bar{x}|}
\end{aligned}$$

$$\geq \limsup_{\substack{t(\bar{x},(\bar{y}_1,\bar{y}_2))\to(\bar{x},(\bar{y}_1,\bar{y}_2))\\ t\uparrow 1}} \frac{-x^*\bar{x}+(y_1^*\bar{y}_1+y_2^*\bar{y}_2)}{2\bar{x}}$$

$$\geq \limsup_{\substack{t(\bar{x},(\bar{y}_1,\bar{y}_2))\to(\bar{x},(\bar{y}_1,\bar{y}_2))\\ t\uparrow 1}} \frac{-x^*\bar{x}-\|(\bar{y}_1,\bar{y}_2)\|_2\|(y_1^*,y_2^*)\|_2}{2\bar{x}}$$

$$= \frac{-x^*-\|(y_1^*,y_2^*)\|_2}{2}$$

$$> 0.$$

This implies that, for any $x^* \in \mathbb{R}$,

$$x^* \in (-\infty, -\|y^*\|_2) \quad \Longrightarrow \quad x^* \notin \widehat{D}^*F(\bar{x}, (\bar{y}_1, \bar{y}_2))(y^*). \tag{3.11}$$

Then, part (b) of (ii) is proved by (3.10) and (3.11) for $\bar{x} > 0$. The case $\bar{x} < 0$ can be similarly proved.

Proof of (iii). Let $x^* \neq 0$. Similarly to the proof of (a) in (ii), by $\bar{x} = 0$ and $\bar{y}_1^2 + \bar{y}_1^2 = 0$, for any $y^* = (y_1^*, y_2^*) \in \mathbb{R}^2$, we calculate

$$\limsup_{\substack{(u,(v_1,v_2))\to(0,(0,0))\\ u\in\mathbb{R} \text{ and } (v_1,v_2)\in F(u)}} \frac{\langle x^*,u-\bar{x}\rangle - \langle (y_1^*,y_2^*),\ (v_1,v_2)-(\bar{y}_1,\bar{y}_2)\rangle}{|u-\bar{x}| + \|(v_1,v_2)-(\bar{y}_1,\bar{y}_2)\|_2}$$

$$= \limsup_{\substack{(u,(v_1,v_2))\to(0,(0,0))\\ u\in\mathbb{R} \text{ and } (v_1,v_2)\in F(u)}} \frac{x^*(u-0)-(y_1^*(v_1-0)+y_2^*(v_2-0))}{|u-\bar{x}| + \|(v_1,v_2)-(\bar{y}_1,\bar{y}_2)\|_2}$$

$$= \limsup_{\substack{(u,(v_1,v_2))\to(0,(0,0))\\ u\in\mathbb{R} \text{ and } (v_1,v_2)\in F(u)}} \frac{x^*u-(y_1^*v_1+y_2^*v_2)}{|u| + \|(v_1,v_2)\|_2}$$

$$\geq \limsup_{\substack{(u,(0,0))\to(0,(0,0))\\ u\in\mathbb{R} \text{ and } (0,0)\in F(u)}} \frac{x^*u-(y_1^*0+y_2^*0)}{|u| + \|(0,0)\|_2}$$

$$= \limsup_{\substack{(u,(0,0))\to(0,(0,0))\\ u\in\mathbb{R} \text{ and } (0,0)\in F(u)}} \frac{x^*u}{|u|}$$

$$= \begin{cases} x^*, & \text{for } t \downarrow 1, \\ -x^*, & \text{for } t \uparrow 1. \end{cases}$$

This implies that

$$\limsup_{\substack{(u,(v_1,v_2))\to(0,(0,0))\\ u\in\mathbb{R} \text{ and } (v_1,v_2)\in F(u)}} \frac{\langle x^*,u-\bar{x}\rangle - \langle (y_1^*,y_2^*),\ (v_1,v_2)-(\bar{y}_1,\bar{y}_2)\rangle}{|u-\bar{x}| + \|(v_1,v_2)-(\bar{y}_1,\bar{y}_2)\|_2} \geq |x^*| > 0.$$

We obtain, for any $y^* = (y_1^*, y_2^*) \in \mathbb{R}^2$,

$$x^* \neq 0 \quad \Longrightarrow \quad x^* \notin \widehat{D}^*F(0, \theta)(y^*). \qquad \square$$

**Example 3.4**. Define $F: (-1,1) \rightrightarrows \mathbb{R}$, by

$$F(x) = \begin{cases} [1,2], & \text{for } -1 < x < 1 \text{ and } x \text{ is rational,} \\ [-2,-1], & \text{for } -1 < x < 1 \text{ and } x \text{ is irrational.} \end{cases}$$

$F$ is a set-valued mapping with closed and convex values. It is clear that $\mathcal{G}(F)^o = \emptyset$. Let $(\bar{x}, \bar{y}) \in \mathcal{G}(F)$. The Mordukhovich derivative $\widehat{D}^*F(\bar{x}, \bar{y})$ of $F$ satisfies that

$$\widehat{D}^*F(\bar{x}, \bar{y})(y^*) \cap (\mathbb{R}\backslash\{0\}) = \emptyset, \text{ for any } y^* \in \mathbb{R}. \tag{3.12}$$

In particular,

(c) if $\bar{x}$ is rational then

$$\widehat{D}^*F(\bar{x}, 2)(y^*) = \{0\}, \text{ for any } y^* \leq 0;$$

$$\widehat{D}^*F(\bar{x}, 1)(y^*) = \{0\}, \text{ for any } y^* \geq 0.$$

(d) if $\bar{x}$ is irrational, then

$$\widehat{D}^*F(\bar{x}, -1)(y^*) = \{0\}, \text{ for any } y^* \leq 0;$$

$$\widehat{D}^*F(\bar{x}, -2)(y^*) = \{0\}, \text{ for any } y^* \geq 0.$$

*Proof*. There are two cases for $\bar{x}$ to be rational or irrational.

Case 1. $(\bar{x}, \bar{y}) \in \mathcal{G}(F)$ with $\bar{x}$ being rational. By definition, $\bar{y} \in [1,2]$. Let $y^* \in \mathbb{R}$ be given. For any $x^* \in \mathbb{R}\backslash\{0\}$, we calculate

$$\limsup_{\substack{(u,v)\to(\bar{x},\bar{y}) \\ u\in(-1,1) \text{ and } v\in F(u)}} \frac{\langle x^*, u-\bar{x}\rangle - \langle y^*, v-\bar{y}\rangle}{|u-\bar{x}| + |v-\bar{y}|}$$

$$\geq \limsup_{\substack{(u,\bar{y})\to(\bar{x},\bar{y}) \\ u \text{ is rational and } \bar{y}\in F(u)}} \frac{\langle x^*, u-\bar{x}\rangle - \langle y^*, \bar{y}-\bar{y}\rangle}{|u-\bar{x}| + |\bar{y}-\bar{y}|}$$

$$= \begin{cases} x^*, & \text{for } u \downarrow \bar{x}, \\ -x^*, & \text{for } u \uparrow \bar{x}. \end{cases}$$

This implies that

$$\limsup_{\substack{(u,v)\to(\bar{x},\bar{y}) \\ u\in(-1,1) \text{ and } v\in F(u)}} \frac{\langle x^*, u-\bar{x}\rangle - \langle y^*, v-\bar{y}\rangle}{|u-\bar{x}| + |v-\bar{y}|} \geq |x^*| > 0.$$

We obtain, for any $y^* \in \mathbb{R}$,

$$x^* \neq 0 \quad \Longrightarrow \quad x^* \notin \widehat{D}^*F(\bar{x}, \bar{y})(y^*). \tag{3.13}$$

Case 2. $(\bar{x}, \bar{y}) \in \mathcal{G}(F)$ with $\bar{x}$ being irrational. By definition, $\bar{y} \in [-2,-1]$. In this case, (3.13) can be similarly proved. This proves part (3.12).

Proof of the first equation of (a). In particular, if $\bar{x}$ is rational and $\bar{y}$ = 2, then, for any $y^* \leq 0$,

$$\limsup_{\substack{(u,v)\to(\bar{x},2)\\ u\in(-1,1) \text{ and } v\in F(u)}} \frac{\langle 0,u-\bar{x}\rangle - \langle y^*,v-2\rangle}{|u-\bar{x}|+|v-2|} = \limsup_{\substack{(u,v)\to(\bar{x},2)\\ u\in(-1,1) \text{ and } v\in F(u)}} \frac{-\langle y^*,v-2\rangle}{|u-\bar{x}|+|v-2|} \leq 0.$$

By (3.12), this proves the first equation of part (a).

Next, we prove the second equation of (a). By the condition that $\bar{x}$ is rational and $\bar{y}$ = 1, when we consider the condition that limit $(u,v) \to (\bar{x},1)$, $u$ must be rational and $v \in F(u) = [1,2]$, as $u$ is near enough to $\bar{x}$. It is because that if $u$ is irrational, then $v \in F(u) = [-2,-1]$ that is impossible for $v$ to approach to 1. Then, for any $y^* \geq 0$, we have

$$\limsup_{\substack{(u,v)\to(\bar{x},1)\\ u\in(-1,1) \text{ and } v\in F(u)}} \frac{\langle 0,u-\bar{x}\rangle - \langle y^*,v-1\rangle}{|u-\bar{x}|+|v-1|} = \limsup_{\substack{(u,v)\to(\bar{x},1)\\ u \text{ is rational and } v\in F(u)}} \frac{-\langle y^*,v-1\rangle}{|u-\bar{x}|+|v-1|} \leq 0.$$

By (3.12), this proves the second equation of part (a). (b) can be similarly proved. □

**Example 3.5.** Mordukhovich Differentiability of the Set-Valued Metric Projection to the closed unit ball in *n*-d Banach Space

Let *n* be a given positive integer. Let $\mathbb{R}^n$ be the ordinary *n*-dimensional real vector space. When $\mathbb{R}^n$ is equipped with the $l_1$-norm, denoted by $\|\cdot\|$, then it becomes the *n*-dimensional $l_1$-Banach space, denoted by $(\mathbb{R}_l^2, \|\cdot\|)$, in which the $l_1$-norm $\|\cdot\|$ on $\mathbb{R}_l^n$ is defined by

$$\|x\| = |x_1| + |x_2| + \ldots + |x_n|, \text{ for any } x = (x_1, x_2, .., x_n) \in \mathbb{R}_l^n.$$

One can prove that the topological dual space of $(\mathbb{R}_l^n, \|\cdot\|)$ is the Banach space $(\mathbb{R}_M^n, \|\cdot\|_M)$, in which the Maximum-norm $\|\cdot\|_M$ on $\mathbb{R}_M^n$ is defined, for any $y^* = (y_1^*, y_2^*, , .., y_n^*,) \in \mathbb{R}_M^n$, by

$$\|y^*\|_M = \max\{|y_1^*|, |y_2^*|, \ldots, |y_n^*|\},$$

As usual, the real pairing between $\mathbb{R}_M^n$ and $\mathbb{R}_l^n$ is written by $\langle\cdot,\cdot\rangle$ and

$$\langle y^*, x\rangle = x_1y_1^* + x_2y_2^* + \cdots + x_ny_n^*, \text{ for any } x = (x_1, x_2, .., x_n) \in \mathbb{R}_l^n \text{ and } y^* = (y_1^*, y_2^*, , .., y_n^*,) \in \mathbb{R}_M^n.$$

$\mathbb{R}_l^n$ and $\mathbb{R}_M^n$ have the same set of elements. In particular, their null element is denoted by $\theta$ = (0, 0, …,0). Let *K* be the positive cone in $\mathbb{R}^n$, which is defined by

$$K = \{x = (x_1, x_2, .., x_n) \in \mathbb{R}^n: x_i \geq 0, \text{for } i = 1, 2, ..., n\}. \tag{3.14}$$

Let $\mathbb{B}$ denote the closed unit ball in $\mathbb{R}_l^n$. Let $P_{\mathbb{B}}: \mathbb{R}_l^n \rightrightarrows \mathbb{B}$ denote the set-valued metric projection operator, which is defined, for any $x \in \mathbb{R}_l^n$, by

$$P_{\mathbb{B}}(x) = \left\{y \in \mathbb{B}: \|x - y\| = \min_{v\in\mathbb{B}}\|x - v\|\right\}.$$

In particular, $P_{\mathbb{B}}(x) = \{x\}$, for any $x \in \mathbb{B}$. Let $x = (x_1, x_2, .., x_n) \in \mathbb{R}_l^n$ with $\|x\| \geq 1$. By Theorem 3.1 in [21], for $y = (y_1, y_2, .., y_n) \in \mathbb{B}$, $y \in P_{\mathbb{B}}(x)$, if and only if $y$ satisfies the following 4 conditions:

(a) $\|y\| = 1$;
(b) $\|x - y\| = \|x\| - 1$;

(c) If $x_i \geq 0$, then $0 \leq y_i \leq x_i \wedge 1$, and if $x_i < 0$, then $x_i \vee (-1) \leq y_i \leq 0$, for $i = 1, 2, \dots, n$.

In particular,

(I) $\frac{x}{\|x\|} \in P_{\mathbb{B}}(x)$, for any $x \in \mathbb{R}_l^n$ with $\|x\| \geq 1$.

(II) If $x \in K$ with $\|x\| \geq 1$, then

$$P_{\mathbb{B}}(x) = \{y \in \mathbb{B}: \textstyle\sum_{i=1}^n y_i = 1, \|x - y\| = \|x\| - 1,\ 0 \leq y_i \leq x_i \wedge 1, \text{ for } 1 \leq i \leq n\}.$$

(III) If $x \in -K$ with $\|x\| \geq 1$, then

$$P_{\mathbb{B}}(x) = \{y \in \mathbb{B}: \textstyle\sum_{i=1}^n y_i = -1, \|x - y\| = \|x\| - 1, x_i \vee (-1) \leq y_i \leq 0, \text{ for } 1 \leq i \leq n\}.$$

Let $\Delta_n$ denote the hyper triangle face in $\mathbb{B}$. More precisely, $\Delta_n$ satisfies that

$$\Delta_n := \{y = (y_1, y_2, \ldots, y_n) \in \mathbb{B} \cap K: 0 \leq y_i \leq 1, \text{for } 1 \leq i \leq n \text{ and } \textstyle\sum_{i=1}^n y_i = 1\}.$$

Let *C* be the convex open subset of *K* that is defined by

$$C := \{x = (x_1, x_2, \ldots, x_n) \in K: x_i > 1, \text{for } 1 \leq i \leq n\}.$$

Then, we have

$$P_{\mathbb{B}}(x) = \Delta_n, \text{ for each } x = (x_1, x_2, \ldots, x_n) \in C.$$

Let arbitrarily given $\bar{x} = (\bar{x}_1, \bar{x}_2, \dots, \bar{x}_n) \in C$. For any $\bar{y} = (\bar{y}_1, \bar{y}_2, \dots, \bar{y}_n) \in P_{\mathbb{B}}(\bar{x}) = \Delta_n$. Let $y^* = (y_1^*, y_2^*, \ldots, y_n^*,) \in \mathbb{R}_M^n \backslash \{\theta\}$. If $\bar{y}_i \neq 1$, for each *i* = 1, 2, …, *n*, then similarly to part (a) in (i) in Theorem 4.2 in [21], the Mordukhovich derivatives of $P_{\mathbb{B}}$ have the following properties:

(I) The equations $y_1^* = y_2^* = \cdots = y_n^*$ do not hold $\Longrightarrow$ $\widehat{D}^* P_{\mathbb{B}_l}(\bar{x}, \bar{y})(y^*) = \emptyset$;

(II) $y_1^* = y_2^* = \cdots = y_n^*$ $\Longrightarrow$ $\widehat{D}^* P_{\mathbb{B}_l}(\bar{x}, \bar{y})(y^*) = \{\theta\}$.

If $\bar{y}_i = 1$, for some *i* = 1, 2, …, *n*, by parts (b−c) in (i) in Theorem 4.2 in [21], we have

$$\widehat{D}^* P_{\mathbb{B}_l}(\bar{x}, (1, 0))(y^*) = \{\theta\} \text{ or } \emptyset.$$

Hence, the results of this example satisfy the results of Theorem 3.1.

## 4. The Closeness and Convexity of Mordukhovich Derivatives

Recall from sections 1 and 3, $(X, \|\cdot\|_X)$ and $(Y, \|\cdot\|_Y)$ are real Banach spaces with topological dual spaces $(X^*, \|\cdot\|_{X^*})$ and $(Y^*, \|\cdot\|_{Y^*})$, respectively. Let $\langle \cdot, \cdot \rangle_X$ denote the real canonical pairing between $X^*$ and *X* and $\langle \cdot, \cdot \rangle_Y$ the real canonical pairing between $Y^*$ and *Y*. Let $\theta_X$ and $\theta_Y$ denote the origins in *X* and *Y*, respectively. Let $\theta_X^*$ and $\theta_Y^*$ denote the null elements in $X^*$ and $Y^*$, respectively. *A* is a nonempty convex open subset of *X*.

**Theorem 4.1**. *Let $F: A \rightrightarrows Y$ be a set-valued mapping. Let $(\bar{x}, \bar{y}) \in \mathcal{G}(F)$ and $y^* \in Y^*$. If $\widehat{D}^* F(\bar{x}, \bar{y})(y^*) \neq \emptyset$, then, $\widehat{D}^* F(\bar{x}, \bar{y})(y^*)$ is a closed (with respect to the $\|\cdot\|_{X^*}$-norm) and convex subset in $X^*$.*

*Proof.* For $(\bar{x}, \bar{y}) \in \mathcal{G}(F)$, let $y^* \in Y^*$ with $\widehat{D}^* F(\bar{x}, \bar{y})(y^*) \neq \emptyset$. Let $x_1^*, x_2^* \in \widehat{D}^* F(\bar{x}, \bar{y})(y^*)$. Then,

$$\limsup_{\substack{(u,v)\to(\bar{x},\bar{y})\\ u\in A \text{ and } v\in F(u)}} \frac{\langle x_i^*, u-\bar{x}\rangle_X - \langle y^*, v-\bar{y}\rangle_Y}{\|u-\bar{x}\|_X + \|v-\bar{y}\|_Y} \leq 0, \text{ for } i = 1, 2. \tag{4.1}$$

For any nonnegative numbers $\alpha_1$and $\alpha_2$ satisfying $\alpha_1 + \alpha_2 = 1$, by (4.1), we estimate

$$\limsup_{\substack{(u,v)\to(\bar{x},\bar{y})\\ u\in A \text{ and } v\in F(u)}} \frac{\langle(\alpha_1 x_1^* + \alpha_2 x_2^*), u-\bar{x}\rangle_X - \langle y^*, v-\bar{y}\rangle_Y}{\|u-\bar{x}\|_X + \|v-\bar{y}\|_Y}$$

$$= \limsup_{\substack{(u,v)\to(\bar{x},\bar{y})\\ u\in A \text{ and } v\in F(u)}} \frac{\langle(\alpha_1 x_1^* + \alpha_2 x_2^*), u-\bar{x}\rangle_X - \langle(\alpha_1 y^* + \alpha_2 y^*), v-\bar{y}\rangle_Y}{\|u-\bar{x}\|_X + \|v-\bar{y}\|_Y}$$

$$= \limsup_{\substack{(u,v)\to(\bar{x},\bar{y})\\ u\in A \text{ and } v\in F(u)}} \frac{\alpha_1(\langle x_1^*, u-\bar{x}\rangle_X - \langle y^*, v-\bar{y}\rangle_Y) + \alpha_2(\langle x_2^*, u-\bar{x}\rangle_X - \langle y^*, v-\bar{y}\rangle_Y)}{\|u-\bar{x}\|_X + \|v-\bar{y}\|_Y}$$

$$\leq \alpha_1 \limsup_{\substack{(u,v)\to(\bar{x},\bar{y})\\ u\in A \text{ and } v\in F(u)}} \frac{\langle x_1^*, u-\bar{x}\rangle_X - \langle y^*, v-\bar{y}\rangle_Y}{\|u-\bar{x}\|_X + \|v-\bar{y}\|_Y} + \alpha_2 \limsup_{\substack{(u,v)\to(\bar{x},\bar{y})\\ u\in A \text{ and } v\in F(u)}} \frac{\langle x_2^*, u-\bar{x}\rangle_X - \langle y^*, v-\bar{y}\rangle_Y}{\|u-\bar{x}\|_X + \|v-\bar{y}\|_Y}$$

$$\leq 0.$$

This implies that $\alpha_1 x_1^* + \alpha_2 x_2^* \in \widehat{D}^* F(\bar{x}, \bar{y})(y^*)$. Next, we prove the $\|\cdot\|_{X^*}$-closeness of $\widehat{D}^* F(\bar{x}, \bar{y})(y^*)$. To this end, suppose that $\{x_n^*\}_{n=1}^{\infty} \subseteq \widehat{D}^* F(\bar{x}, \bar{y})(y^*)$ and $x_n^* \to x^*$, as $n \to \infty$ in $X^*$, with respect to $\|\cdot\|_{X^*}$-norm. For any $\varepsilon > 0$, there is a positive integer $N$ such that

$$\|x^* - x_n^*\|_{X^*} < \varepsilon, \text{ for any } n > N. \tag{4.2}$$

Then, for any given $n > N$, by (4.2), we estimate

$$\limsup_{\substack{(u,v)\to(\bar{x},\bar{y})\\ u\in A \text{ and } v\in F(u)}} \frac{\langle x^*, u-\bar{x}\rangle_X - \langle y^*, v-\bar{y}\rangle_Y}{\|u-\bar{x}\|_X + \|v-\bar{y}\|_Y}$$

$$= \limsup_{\substack{(u,v)\to(\bar{x},\bar{y})\\ u\in A \text{ and } v\in F(u)}} \frac{\langle x^* - x_n^*, u-\bar{x}\rangle_X + \langle x_n^*, u-\bar{x}\rangle_X - \langle y^*, v-\bar{y}\rangle_Y}{\|u-\bar{x}\|_X + \|v-\bar{y}\|_Y}$$

$$\leq \limsup_{\substack{(u,v)\to(\bar{x},\bar{y})\\ u\in A \text{ and } v\in F(u)}} \frac{\langle x^* - x_n^*, u-\bar{x}\rangle_X}{\|u-\bar{x}\|_X + \|v-\bar{y}\|_Y} + \limsup_{\substack{(u,v)\to(\bar{x},\bar{y})\\ u\in A \text{ and } v\in F(u)}} \frac{\langle x_n^*, u-\bar{x}\rangle_X - \langle y^*, v-\bar{y}\rangle_Y}{\|u-\bar{x}\|_X + \|v-\bar{y}\|_Y}$$

$$\leq \limsup_{\substack{(u,v)\to(\bar{x},\bar{y})\\ u\in A \text{ and } v\in F(u)}} \frac{\langle x^* - x_n^*, u-\bar{x}\rangle_X}{\|u-\bar{x}\|_X + \|v-\bar{y}\|_Y} + 0$$

$$\leq \limsup_{\substack{(u,v)\to(\bar{x},\bar{y})\\ u\in A \text{ and } v\in F(u)}} \frac{\|x^* - x_n^*\|_{X^*} \|u-\bar{x}\|_X}{\|u-\bar{x}\|_X + \|v-\bar{y}\|_Y}$$

$$\leq \limsup_{\substack{(u,v)\to(\bar{x},\bar{y})\\ u\in A \text{ and } v\in F(u)}} \frac{\|x^* - x_n^*\|_{X^*} \|u-\bar{x}\|_X}{\|u-\bar{x}\|_X}$$

$= \|x^* - x_n^*\|_{X^*}$

$< \varepsilon.$

This implies that

$$\limsup_{\substack{(u,v)\to(\bar{x},\bar{y})\\ u\in A \text{ and } v\in F(u)}} \frac{\langle x^*, u-\bar{x}\rangle_X - \langle y^*, v-\bar{y}\rangle_Y}{\|u-\bar{x}\|_X + \|v-\bar{y}\|_Y} \leq 0.$$

Hence $x^* \in \widehat{D}^*F(\bar{x},\bar{y})(y^*)$; and therefore, $\widehat{D}^*F(\bar{x},\bar{y})(y^*)$ is $\|\cdot\|_{X^*}$-closed in $X^*$. □

Next, we consider the Mordukhovich differentiability of the standard metric projection operator in general Banach spaces. In Example 3.5, and in Proposition 3.1 and Theorem 4.1 in [21], we find that the metric projection operator to the closed unit ball is indeed a set-valued mapping. However, in the following proposition, we will find that when the metric projection operator $P_K$ projects to the positive cone $K$ in $\mathbb{R}_l^n$ defined by (3.14), $P_K$ is a single-valued mapping.

Let $n$ be a given positive integer. Let $(\mathbb{R}_l^n, \|\cdot\|)$ be the $n$-dimensional $l_1$-Banach space with topological dual space $(\mathbb{R}_M^n, \|\cdot\|_M)$ defined in Example 3.5. Let $K$ be the positive cone in $\mathbb{R}_l^n$ defined in (3.14). Let $P_K$ denote the metric projection operator, which is defined, for any $x = (x_1, x_2, .., x_n) \in \mathbb{R}_l^n$, by

$$P_K(x) = \left\{y \in K: \|x-y\| = \min_{v\in K}\|x-v\|\right\}. \tag{4.3}$$

The topological interior $K^o$ and the boundary $\partial K$ satisfy the following equations.

$$K^o = \{y = (y_1, y_2, .., y_n) \in K: y_1 y_2 \dots y_n > 0\} \text{ and } \partial K = \{y = (y_1, y_2, .., y_n) \in K: y_1 y_2 \dots y_n = 0\}.$$

Let $\preccurlyeq_K$ denote the partial order on $\mathbb{R}_l^n$ generated by the positive cone $K$, which is defined, for any $x = (x_1, x_2, .., x_n) \in \mathbb{R}_l^n$ and $u = (u_1, u_2, .., u_n) \in \mathbb{R}_l^n$, by

$$x \preccurlyeq_K u \quad \text{if and only if} \quad x_i \leq u_i, \text{ for each } i = 1, 2, \dots, n.$$

**Proposition 4.2**. *The metric projection operator* $P_K: \mathbb{R}_l^n \to K$ *defined in* (3.14) *is a single-valued mapping and, for any* $x = (x_1, x_2, .., x_n) \in \mathbb{R}_l^n$, $P_K(x) = ((P_K(x))_1, (P_K(x))_2, \dots (P_K(x))_n) \in K$ *has the following representation*

$$(P_K(x))_i = \begin{cases} x_i, & \text{if } x_i > 0, \\ 0, & \text{if } x_i \leq 0, \end{cases} \quad \text{for each } i = 1, 2, \dots, n. \tag{4.4}$$

*Proof*. For any $v = (v_1, v_2, .., v_n) \in K$ satisfying $v_i \geq 0$, for $i = 1, 2, \dots, n$, we calculate

$$\|x - v\| = \sum_{i=1}^n |x_i - v_i|$$

$$= \sum_{x_i>0} |x_i - v_i| + \sum_{x_i\leq 0} |x_i - v_i|$$

$$\geq \sum_{x_i>0} |x_i - x_i| + \sum_{x_i\leq 0} |x_i - 0|.$$

By definition, this implies (4.4). □

In the following proposition, we will show that at some points in Banach space $\mathbb{R}_l^n$, the metric projection $P_K$ is not Fréchet differentiable**.** Meanwhile, the Mordukhovich derivative is a closed and convex subset of $\mathbb{R}^n_M$, which demonstrates the results of Theorem 4.1 in this paper.

**Proposition 4.3**. *The metric projection operator* $P_K\colon \mathbb{R}^n_l \to K$ *defined in* (3.14) *has the following Fréchet and Mordukhovich differentiability*.

(i) $P_K$ *is not Fréchet differentiable at any point* $\bar{x} \in \partial K$.
(ii) *Let* $y^* \in K$. *The Mordukhovich derivative satisfies the following equation.*

$$\widehat{D}^* P_K(\theta, \theta)(y^*) = \{x^* \in \mathbb{R}^n_M : \theta \preccurlyeq_K x^* \preccurlyeq_K y^*\}.$$

Notice that $\{x^* \in \mathbb{R}^n_M : \theta \preccurlyeq_K x^* \preccurlyeq_K y^*\}$ is indeed a nonempty closed and convex subset of $\mathbb{R}^n_M$.

*Proof*. Proof of (i). Let $\bar{x} = (\bar{x}_1, \bar{x}_2, \ldots, \bar{x}_n) \in \partial K$ with $\bar{x}_1 \bar{x}_2 \ldots \bar{x}_n = 0$. The proof is divided to two cases.

Case 1. Suppose that there are positive integers *k* and *m* such that $\bar{x}_k = 0$ and $\bar{x}_m > 0$. For each *j* = 1, 2, …, *n*, let $\beta_j$ denote the vector in $\mathbb{R}^n_l$, in which the $j^{th}$ coordinate is 1 and all other coordinates are 0. For real number *t*, define $u(t) = t\beta_k + t\beta_m$. Assume, by the way of contradiction, that $\nabla P_K(\bar{x})$ exists, which is a continuous and linear mapping from $\mathbb{R}^n_l$ to $\mathbb{R}^n_l$. Then, we calculate

$$\begin{aligned}
0 &= \lim_{u \to 0} \frac{P_K(\bar{x}+u) - P_K(\bar{x}) - \nabla P_K(\bar{x})(u)}{\|u\|} \\
&= \lim_{t \downarrow 0} \frac{P_K(\bar{x}+u(\mathrm{t})) - P_K(\bar{x}) - \nabla P_K(\bar{x})(u(t))}{\|u(t)\|} \\
&= \lim_{t \downarrow 0} \frac{t(\beta_k + \beta_m) - t\nabla P_K(\bar{x})(\beta_k + \beta_m)}{2t} \\
&= \frac{(\beta_k + \beta_m) - \nabla P_K(\bar{x})(\beta_k + \beta_m)}{2}
\end{aligned}$$

This implies

$$\nabla P_K(\bar{x})(\beta_k + \beta_m) = \beta_k + \beta_m. \tag{4.5}$$

On the other hand, we have

$$\begin{aligned}
0 &= \lim_{u \to 0} \frac{P_K(\bar{x}+u) - P_K(\bar{x}) - \nabla P_K(\bar{x})(u)}{\|u\|} \\
&= \lim_{t \uparrow 0} \frac{P_K(\bar{x}+u(\mathrm{t})) - P_K(\bar{x}) - \nabla P_K(\bar{x})(u(t))}{\|u(t)\|} \\
&= \lim_{t \uparrow 0, |t| < \bar{x}_m} \frac{t(\theta + \beta_m) - t\nabla P_K(\bar{x})(\beta_k + \beta_m)}{2|t|} \\
&= \frac{-\beta_m + \nabla P_K(\bar{x})(\beta_k + \beta_m)}{2}.
\end{aligned}$$

This implies

$$\nabla P_K(\bar{x})(\beta_k + \beta_m) = \beta_m. \tag{4.6}$$

By $\beta_k \neq \theta$, (4.5) and (4.6) imply that, for $\bar{x} = (\bar{x}_1, \bar{x}_2, \ldots, \bar{x}_n) \in \partial K$, $P_K$ is not Fréchet differentiable at $\bar{x}$.

Case 2. Let $\bar{x} = \theta$. At first, similarly to the proof of (4.5), we have

$$0 = \lim_{u \to 0} \frac{P_K(\theta+u) - P_K(\theta) - \nabla P_K(\theta)(u)}{\|u\|}$$

$$= \lim_{t \downarrow 0} \frac{P_K(\theta+u(\mathrm{t})) - P_K(\theta) - \nabla P_K(\theta)(u(t))}{\|u(t)\|}$$

$$= \lim_{t \downarrow 0} \frac{t(\beta_k+\beta_m) - t\nabla P_K(\theta)(\beta_k+\beta_m)}{2t}$$

$$= \frac{(\beta_k+\beta_m) - \nabla P_K(\theta)(\beta_k+\beta_m)}{2}$$

This implies

$$\nabla P_K(\theta)(\beta_k + \beta_m) = \beta_k + \beta_m. \tag{4.7}$$

Similarly to the proof of (4.6), we have

$$0 = \lim_{u \to 0} \frac{P_K(\theta+u) - P_K(\theta) - \nabla P_K(\theta)(u)}{\|u\|}$$

$$= \lim_{t \uparrow 0} \frac{P_K(\theta+u(\mathrm{t})) - P_K(\theta) - \nabla P_K(\theta)(u(t))}{\|u(t)\|}$$

$$= \lim_{t \uparrow 0, |t| < \bar{x}_m} \frac{t(\theta) - t\nabla P_K(\theta)(\beta_k+\beta_m)}{2|t|}$$

$$= \frac{\theta + \nabla P_K(\theta)(\beta_k+\beta_m)}{2}.$$

This implies

$$\nabla P_K(\bar{x})(\beta_k + \beta_m) = \theta. \tag{4.8}$$

By $\beta_k + \beta_m \neq \theta$, (4.7) and (4.8) imply that, for $\bar{x} = \theta \in \partial K$, $P_K$ is not Fréchet differentiable at $\theta$.

Proof of (ii). Let $y^* = (y_1^*, y_2^*, \ldots, y_n^*,) \in K$ and $x^* = (x_1^*, x_2^*, \ldots, x_n^*,) \in K$. Suppose that $\theta \preccurlyeq_K x^* \preccurlyeq_K y^*$. Let $\bar{x} = \theta$. Let $u = (u_1, u_2, \ldots, u_n) \in \mathbb{R}_l^n$ and $v = (v_1, v_2, \ldots, v_n) = P_K(u)$. By the representation of $P_K$, it implies that if $u_i \leq 0$, then $v_i = P_K(u)_i = 0$ and if $u_i > 0$, then $v_i = P_K(u)_i = u_i$. By the assumption that $0 \leq x_i^* \leq y_i^*$, for each $i$ = 1, 2, …, $n$, we estimate

$$\limsup_{\substack{(u,v) \to (\bar{x},\bar{x}) \\ v = P_K(u)}} \frac{\langle x^*, u - \bar{x} \rangle - \langle y^*, v - \bar{x} \rangle}{\|u - \bar{x}\| + \|v - \bar{x}\|_M}$$

$$= \limsup_{\substack{(u,v) \to (\theta,\theta) \\ v = P_K(u)}} \frac{\langle x^*, u - \theta \rangle - \langle y^*, v - \theta \rangle}{\|u - \theta\| + \|v - \theta\|_M}$$

$$= \limsup_{\substack{(u,v) \to (\theta,\theta) \\ v = P_K(u)}} \frac{\langle x^*, u \rangle - \langle y^*, v \rangle}{\|u - \theta\| + \|v - \theta\|_M}$$

$$= \limsup_{\substack{(u,v) \to (\theta,\theta) \\ v = P_K(u)}} \frac{\sum_{u_i > 0}(x_i^* u_i - y_i^* v_i) + \sum_{u_i \leq 0}(x_i^* u_i - y_i^* v_i)}{\|u - \theta\| + \|v - \theta\|_M}$$

$$= \limsup_{\substack{(u,v)\to(\theta,\theta)\\ v=P_K(u)}} \frac{\sum_{u_i>0}(x_i^* u_i - y_i^* u_i) + \sum_{u_i\le 0}(x_i^* u_i - y_i^* 0)}{\|u-\theta\| + \|v-\theta\|_M}$$

$$= \limsup_{\substack{(u,v)\to(\theta,\theta)\\ v=P_K(u)}} \frac{\sum_{u_i>0}(x_i^* - y_i^*)u_i + \sum_{u_i\le 0}(x_i^* u_i)}{\|u-\theta\| + \|v-\theta\|_M}$$

$$\le 0.$$

This implies that, for given $y^* = (y_1^*, y_2^*, \dots, y_n^*,) \in K$,

$$x^* \in \widehat{D}^* P_K(\theta,\theta)(y^*), \text{ for any } \theta \preccurlyeq_K x^* \preccurlyeq_K y^*. \qquad (4.9)$$

Suppose $\theta \not\preccurlyeq_K x^*$. Then, there is a positive number *m* with $x_m^* < 0$. By the proof of (4.9), we have

$$\limsup_{\substack{(u,v)\to(\bar{x},\bar{x})\\ v=P_K(u)}} \frac{\langle x^*, u-\bar{x}\rangle - \langle y^*, v-\bar{x}\rangle}{\|u-\bar{x}\| + \|v-\bar{x}\|_M}$$

$$= \limsup_{\substack{(u,v)\to(\theta,\theta)\\ v=P_K(u)}} \frac{\sum_{u_i>0}(x_i^* - y_i^*)u_i + \sum_{u_i\le 0}(x_i^* u_i)}{\|u-\theta\| + \|v-\theta\|_M}$$

$$\ge \limsup_{\substack{(t\beta_m,\theta)\to(\theta,\theta)\\ \theta=P_K(t\beta_m),\, t\uparrow 0}} \frac{\sum_{u_i>0}(x_i^* - y_i^*)0 + x_m^* t}{|t| + 0}$$

$$\ge \limsup_{\substack{(t\beta_m,\theta)\to(\theta,\theta)\\ \theta=P_K(t\beta_m),\, t\uparrow 0}} \frac{-x_m^*}{1}$$

$$= -x_m^* > 0.$$

This implies that for given $y^* = (y_1^*, y_2^*, \dots, y_n^*,) \in K$,

$$x^* \notin \widehat{D}^* P_K(\theta,\theta)(y^*), \text{ for any } \theta \not\preccurlyeq_K x^*. \qquad (4.10)$$

Next, suppose $x^* \not\preccurlyeq_K y^*$. Then, there is a positive number *k* with $x_k^* > y_k^*$. By the proof of (4.10), we have

$$\limsup_{\substack{(u,v)\to(\bar{x},\bar{x})\\ v=P_K(u)}} \frac{\langle x^*, u-\bar{x}\rangle - \langle y^*, v-\bar{x}\rangle}{\|u-\bar{x}\| + \|v-\bar{x}\|_M}$$

$$= \limsup_{\substack{(u,v)\to(\theta,\theta)\\ v=P_K(u)}} \frac{\sum_{u_i>0}(x_i^* - y_i^*)u_i + \sum_{u_i\le 0}(x_i^* u_i)}{\|u-\theta\| + \|v-\theta\|_M}$$

$$\ge \limsup_{\substack{(t\beta_k,t\beta_k)\to(\theta,\theta)\\ t\beta_k=P_K(t\beta_m),\, t\downarrow 0}} \frac{(x_k^* - y_k^*)t + \sum_{u_i\le 0}(x_i^* 0)}{|t| + |t|}$$

$$\geq \frac{x_k^* - y_k^*}{2}$$

$$> 0.$$

This implies that for given $y^* = (y_1^*, y_2^*, \ldots, y_n^*,) \in K$,

$$x^* \notin \widehat{D}^* P_K(\theta, \theta)(y^*), \text{ for any } x^* \not\preccurlyeq_K y^*. \tag{4.11}$$

This theorem is proved by (4.9), (4.10) and (4.11). □

**Notations 4.4**. Let $F: A \rightrightarrows Y$ be a set-valued mapping. Let $(\bar{x}, \bar{y}) \in \mathcal{G}(F)$. We write

$$\widehat{D}^* F(\bar{x}, \bar{y}) = \{y^* \in Y^* : \widehat{D}^* F(\bar{x}, \bar{y})(y^*) \neq \emptyset\}.$$

It is clear that, for any $(\bar{x}, \bar{y}) \in \mathcal{G}(F)$, we have $\theta_{X^*} \in \widehat{D}^* F(\bar{x}, \bar{y})(\theta_{Y^*})$, which implies that

$$\theta_{Y^*} \in \widehat{D}^* F(\bar{x}, \bar{y}), \text{ for any } (\bar{x}, \bar{y}) \in \mathcal{G}(F).$$

For given $(\bar{x}, \bar{y}) \in \mathcal{G}(F)$, if $\widehat{D}^* F(\bar{x}, \bar{y}) \supsetneqq \{\theta_{Y^*}\}$, that is, there is $y^* \in Y^* \backslash \{\theta_{Y^*}\}$ such that $\widehat{D}^* F(\bar{x}, \bar{y})(y^*) \neq \emptyset$, then $F$ is said to be Mordukhovich differentiable at $(\bar{x}, \bar{y}) \in \mathcal{G}(F)$.

**Proposition 4.5**. *Let $F: A \rightrightarrows Y$ be a set-valued mapping. Let $(\bar{x}, \bar{y}) \in \mathcal{G}(F)$. If F is Mordukhovich differentiable at $(\bar{x}, \bar{y})$, then, $\widehat{D}^* F(\bar{x}, \bar{y})$ is a convex cone in $Y^*$ with vertex at $\theta_{Y^*}$.*

*Proof*. Let $y^*, z^* \in Y^* \backslash \{\theta_{Y^*}\}$ such that $\widehat{D}^* F(\bar{x}, \bar{y})(y^*) \neq \emptyset$ and $\widehat{D}^* F(\bar{x}, \bar{y})(z^*) \neq \emptyset$. Then, there are $x^*, w^* \in X^*$ such that

$$x^* \in \widehat{D}^* F(\bar{x}, \bar{y})(y^*) \quad \text{and} \quad w^* \in \widehat{D}^* F(\bar{x}, \bar{y})(z^*).$$

By the condition, these are equivalent to

$$\limsup_{\substack{(u,v) \to (\bar{x}, \bar{y}) \\ u \in A \text{ and } v \in F(u)}} \frac{\langle x^*, u - \bar{x} \rangle_X - \langle y^*, v - \bar{y} \rangle_Y}{\|u - \bar{x}\|_X + \|v - \bar{y}\|_Y} \leq 0 \quad \text{and} \quad \limsup_{\substack{(u,v) \to (\bar{x}, \bar{y}) \\ u \in A \text{ and } v \in F(u)}} \frac{\langle w^*, u - \bar{x} \rangle_X - \langle z^*, v - \bar{y} \rangle_Y}{\|u - \bar{x}\|_X + \|v - \bar{y}\|_Y} \leq 0.$$

Let $\alpha, \beta$ be nonnegative numbers with $\alpha + \beta = 1$. We have

$$\limsup_{\substack{(u,v) \to (\bar{x}, \bar{y}) \\ u \in A \text{ and } v \in F(u)}} \frac{\langle \alpha x^* + \beta w^*, u - \bar{x} \rangle_X - \langle \alpha y^* + \beta z^*, v - \bar{y} \rangle_Y}{\|u - \bar{x}\|_X + \|v - \bar{y}\|_Y}$$

$$\leq \alpha \limsup_{\substack{(u,v) \to (\bar{x}, \bar{y}) \\ u \in A \text{ and } v \in F(u)}} \frac{\langle x^*, u - \bar{x} \rangle_X - \langle y^*, v - \bar{y} \rangle_Y}{\|u - \bar{x}\|_X + \|v - \bar{y}\|_Y} + \beta \limsup_{\substack{(u,v) \to (\bar{x}, \bar{y}) \\ u \in A \text{ and } v \in F(u)}} \frac{\langle w^*, u - \bar{x} \rangle_X - \langle z^*, v - \bar{y} \rangle_Y}{\|u - \bar{x}\|_X + \|v - \bar{y}\|_Y}$$

$$\leq 0.$$

This implies that $\alpha x^* + \beta w^* \in \widehat{D}^* F(\bar{x}, \bar{y})(\alpha y^* + \beta z^*)$; and therefore,

$$\alpha y^* + \beta z^* \in \widehat{D}^* F(\bar{x}, \bar{y}).$$

This proves the convexity of $\widehat{D}^* F(\bar{x}, \bar{y})$. Then, we prove that $\widehat{D}^* F(\bar{x}, \bar{y})$ is a cone in $Y^*$. To this end, let

$y^* \in Y^* \backslash \{\theta_{Y^*}\}$ and $x^* \in X^*$ such that

$$x^* \in \widehat{D}^* F(\bar{x}, \bar{y})(y^*), \text{ that is } \limsup_{\substack{(u,v)\to(\bar{x},\bar{y}) \\ u\in A \text{ and } v\in F(u)}} \frac{\langle x^*, u-\bar{x}\rangle_X - \langle y^*, v-\bar{y}\rangle_Y}{\|u-\bar{x}\|_X + \|v-\bar{y}\|_Y} \leq 0.$$

Let $\gamma$ be a nonnegative number. We have

$$\limsup_{\substack{(u,v)\to(\bar{x},\bar{y}) \\ u\in A \text{ and } v\in F(u)}} \frac{\langle \gamma x^*, u-\bar{x}\rangle_X - \langle \gamma y^*, v-\bar{y}\rangle_Y}{\|u-\bar{x}\|_X + \|v-\bar{y}\|_Y}$$

$$= \gamma \limsup_{\substack{(u,v)\to(\bar{x},\bar{y}) \\ u\in A \text{ and } v\in F(u)}} \frac{\langle x^*, u-\bar{x}\rangle_X - \langle y^*, v-\bar{y}\rangle_Y}{\|u-\bar{x}\|_X + \|v-\bar{y}\|_Y}$$

$$\leq 0.$$

This proves that $\widehat{D}^* F(\bar{x}, \bar{y})$ is a cone in $Y^*$. □

**Proposition 4.6.** *Let* $P_K: \mathbb{R}^n_l \to K$ *be the metric projection operator defined in* (3.14). *Then,*

$$\widehat{D}^* P_K(\theta, \theta) = K. \tag{4.12}$$

*Proof*. By part (ii) in Proposition 4.3, for any $y^* \in K$, we have

$$\widehat{D}^* P_K(\theta, \theta)(y^*) = \{x^* \in \mathbb{R}^n_M : \theta \preccurlyeq_K x^* \preccurlyeq_K y^*\} \coloneqq [\theta, y^*]_{\preccurlyeq_K}.$$

This implies that, for any $y^* \in K$ we have $y^* \in \widehat{D}^* P_K(\theta, \theta)$. Hence

$$K \subseteq \widehat{D}^* P_K(\theta, \theta). \tag{4.13}$$

If $y^* = (y_1^*, y_2^*, \ldots, y_n^*,) \notin K$, then there is *m* such that $y_m^* < 0$. We calculate

$$\limsup_{\substack{(u,v)\to(\bar{x},\bar{x}) \\ v=P_K(u)}} \frac{\langle x^*, u-\bar{x}\rangle - \langle y^*, v-\bar{x}\rangle}{\|u-\bar{x}\| + \|v-\bar{x}\|_M}$$

$$= \limsup_{\substack{(u,v)\to(\theta,\theta) \\ v=P_K(u)}} \frac{\langle x^*, u-\theta\rangle - \langle y^*, v-\theta\rangle}{\|u-\theta\| + \|v-\theta\|_M}$$

$$= \limsup_{\substack{(u,v)\to(\theta,\theta) \\ v=P_K(u)}} \frac{\langle x^*, u\rangle - \langle y^*, v\rangle}{\|u-\theta\| + \|v-\theta\|_M}$$

$$= \limsup_{\substack{(u,v)\to(\theta,\theta) \\ v=P_K(u)}} \frac{\sum_{u_i>0}(x_i^* u_i - y_i^* v_i) + \sum_{u_i\leq 0}(x_i^* u_i - y_i^* v_i)}{\|u-\theta\| + \|v-\theta\|_M}$$

$$= \limsup_{\substack{(u,v)\to(\theta,\theta) \\ v=P_K(u)}} \frac{\sum_{u_i>0}(x_i^* u_i - y_i^* u_i) + \sum_{u_i\leq 0}(x_i^* u_i - y_i^* 0)}{\|u-\theta\| + \|v-\theta\|_M}$$

$$= \limsup_{\substack{(u,v)\to(\theta,\theta)\\ v=P_K(u)}} \frac{\sum_{u_i>0}(x_i^*-y_i^*)u_i+\sum_{u_i\le 0}(x_i^*u_i)}{\|u-\theta\| + \|v-\theta\|_M}. \qquad (4.14)$$

Next, let arbitrarily $x^* = (x_1^*, x_2^*, .., x_n^*,) \in \mathbb{R}_M^n$. We consider the following two cases.

Case 1. $x^* \notin K$. There is a positive integer $k$ with $x_k^* < 0$. For $t < 0$, let $u(t) = (u_1, u_2, .., u_n) \in \mathbb{R}_l^n$ by

$$u_i = \begin{cases} t, & \text{if } i = k, \\ 0, & \text{if } i \neq k, \end{cases} \text{ for each } i = 1, 2, \ldots, n.$$

Substituting $u(t)$ in (4.14), we have

$$\limsup_{\substack{(u,v)\to(\theta,\theta)\\ v=P_K(u)}} \frac{\langle x^*, u-\theta\rangle - \langle y^*, v-\theta\rangle}{\|u-\theta\| + \|v-\theta\|_M}$$

$$= \limsup_{\substack{(u,v)\to(\theta,\theta)\\ v=P_K(u)}} \frac{\sum_{u_i>0}(x_i^*u_i-y_i^*v_i)+\sum_{u_i\le 0}(x_i^*u_i-y_i^*v_i)}{\|u-\theta\| + \|v-\theta\|_M}$$

$$\ge \limsup_{\substack{u(t)\to\theta,t\uparrow 0\\ v=P_K(u(t))=\theta}} \frac{x_k^*t}{|t|}$$

$$= -x_k^* > 0.$$

This implies that, for $x^* = (x_1^*, x_2^*, .., x_n^*,) \neq \theta$,

$$x^* \notin K \quad \Longrightarrow \quad x^* \notin \widehat{D}^*P_K(\theta,\theta)(y^*). \qquad (4.15)$$

Case 2. $x^* = (x_1^*, x_2^*, .., x_n^*,) \in K$, that is $x_i^* \geq 0$, for each $i = 1, 2, \ldots, n$. In this case, for $t > 0$, let $u(t) = (u_1, u_2, .., u_n) \in \mathbb{R}_l^n$ by

$$u_i = \begin{cases} t, & \text{if } i = m, \\ 0, & \text{if } i \neq m, \end{cases} \text{ for each } i = 1, 2, \ldots, n.$$

Substituting $u(t)$ in (4.14), we have

$$\limsup_{\substack{(u,v)\to(\theta,\theta)\\ v=P_K(u)}} \frac{\langle x^*, u-\theta\rangle - \langle y^*, v-\theta\rangle}{\|u-\theta\| + \|v-\theta\|_M}$$

$$= \limsup_{\substack{(u,v)\to(\theta,\theta)\\ v=P_K(u)}} \frac{\sum_{u_i>0}(0-y_i^*)u_i+\sum_{u_i\le 0}(0)}{\|u-\theta\| + \|v-\theta\|_M}$$

$$\ge \limsup_{\substack{u(t)\to\theta,t\downarrow 0\\ v=P_K(u(t))=u(t)}} \frac{(x_m^*-y_m^*)t}{2t}$$

$$= \frac{x_m^*-y_m^*}{2} > 0 \ (x_m^* \geq 0 \text{ and } y_m^* < 0)$$

This implies that

$$x^* \in K \quad \Longrightarrow \quad x^* \notin \widehat{D}^* P_K(\theta, \theta)(y^*). \tag{4.16}$$

By (4.15) and (4.16), we obtain that

$$y^* \notin K \quad \Longrightarrow \quad \widehat{D}^* P_K(\theta, \theta)(y^*) = \emptyset. \tag{4.17}$$

By definition, (4.12) is proved by (4.13) and (4.17).

## 5. Applications to Locally Variational Inequalities in Banach Spaces

In this section, we first introduce the concept of locally variational inequality problems. Then, we investigate the connection between the solutions of some locally variational inequality problems and Mordukhovich derivatives with respect to the same set-valued mapping. We start at the following definition of locally variational inequality problems.

**Definition 5.1**. Let $F: A \rightrightarrows Y$ be a set-valued mapping. Let $\bar{x} \in A$ and $\bar{y} \in F(\bar{x})$. Let $C$ be an open neighborhood of $\bar{x}$ in $X$, with $C \subseteq A$. The locally variational inequality problem associated with $F$ and $C$, at point $(\bar{x}, \bar{y}) \in \mathcal{G}(F)$ is to find $x^* \in X^*$ and $y^* \in Y^*$ such that

$$\langle y^*, v - \bar{y} \rangle_Y \geq \langle x^*, u - \bar{x} \rangle_X, \text{ for any } u \in C \text{ and } v \in F(u). \tag{5.1}$$

Such a locally variational inequality problem defined by (5.1) is denoted by LVI($F(\bar{x}, \bar{y})$, $C$). The point $(x^*, y^*) \in X^* \times Y^*$ satisfying (5.1) is called a solution to the problem LVI($F(\bar{x}, \bar{y})$, $C$). The collection of all solutions to the problem LVI($F(\bar{x}, \bar{y})$, $C$) is denoted by SLVI($F(\bar{x}, \bar{y})$, $C$), which has the following properties.

**Lemma 5.2**. *Let $F: A \rightrightarrows Y$ be a set-valued mapping. Let $\bar{x} \in A$ and $\bar{y} \in F(\bar{x})$. Let $C \subseteq A$ be an open neighborhood of $\bar{x}$ in $X$.*

(i) *Let $y^* \in Y^*$. Define a subset of $X^*$ by*

$$\{x^* \in X^*: (x^*, y^*) \in \text{SLVI}(F(\bar{x}, \bar{y}), C)\}. \tag{5.2}$$

*Then the above set given by* (5.2) *is either empty or a nonempty convex subset of $X^*$. Furthermore, if $C$ is $\|\cdot\|_X$-bounded, then the set defined by* (5.2) *is $\|\cdot\|_{X^*}$-closed*;

(ii) *Let $x^* \in X^*$. Define a subset of $Y^*$ by*

$$\{y^* \in Y^*: (x^*, y^*) \in \text{SLVI}(F(\bar{x}, \bar{y}), C)\}. \tag{5.3}$$

*Then the above set defined by* (5.3) *is either empty or a nonempty convex subset of $Y^*$*;

(iii) *Let $(x^*, y^*), (u^*, v^*) \in \text{SLVI}(F(\bar{x}, \bar{y}), C)$. Let $\alpha, \beta$ be nonnegative with $\alpha + \beta = 1$. Then,*

$$(\alpha x^* + \beta u^*, \alpha y^* + \beta v^*) \in \text{SLVI}(F(\bar{x}, \bar{y}), C).$$

(iv) *Let $(x^*, y^*) \in \text{SLVI}(F(\bar{x}, \bar{y}), C)$. Then,*

$$(\lambda x^*, \lambda y^*) \in \text{SLVI}(F(\bar{x}, \bar{y}), C), \textit{ for any } \lambda \in \mathbb{R}_+.$$

*Proof*. The proof of this lemma is similar to the proof of Theorem 4.1 and it is omitted here. □

**Theorem 5.3**. *Let $F: A \rightrightarrows Y$ be a set-valued mapping. Let $\bar{x} \in A$ and $\bar{y} \in F(\bar{x})$. Let $C \subseteq A$ be an open neighborhood of $\bar{x}$ in X. Let $(x^*, y^*) \in X^* \times Y^*$. If $(x^*, y^*)$ is a solution to the problem* LVI($F(\bar{x}, \bar{y})$, $C$), *then $x^*$ and $y^*$ satisfy*

$$x^* \in \widehat{D}^*F(\bar{x}, \bar{y})(y^*).$$

*Proof*. We calculate

$$\limsup_{\substack{(u,v)\to(\bar{x},\bar{y}) \\ u\in A \text{ and } v\in F(u)}} \frac{\langle x^*, u-\bar{x}\rangle_X - \langle y^*, v-\bar{y}\rangle_Y}{\|u-\bar{x}\|_X + \|v-\bar{y}\|_Y}$$

$$= \limsup_{\substack{(u,v)\to(\bar{x},\bar{y}) \\ u\in U\cap C \text{ and } v\in F(u)\cap V}} \frac{\langle x^*, u-\bar{x}\rangle_X - \langle y^*, v-\bar{y}\rangle_Y}{\|u-\bar{x}\|_X + \|v-\bar{y}\|_Y}$$

$$\leq \limsup_{\substack{(u,v)\to(\bar{x},\bar{y}) \\ u\in U\cap C \text{ and } v\in F(u)\cap V}} \frac{0}{\|u-\bar{x}\|_X + \|v-\bar{y}\|_Y}$$

$$= 0.$$

This implies $x^* \in \widehat{D}^*F(\bar{x}, \bar{y})(y^*)$. □

The following consequence of Theorem 5.3 is considered as an application of Mordukhovich derivatives to locally variational inequality problems.

**Corollary 5.4**. *Let $F: A \rightrightarrows Y$ be a set-valued mapping. Let $\bar{x} \in A$ and $\bar{y} \in F(\bar{x})$. Let $(x^*, y^*) \in X^* \times Y^*$. Then, for any open neighborhood C of $\bar{x}$ in X, with $C \subseteq A$, we have*

$$x^* \notin \widehat{D}^*F(\bar{x}, \bar{y})(y^*) \quad \Longrightarrow \quad (x^*, y^*) \notin \text{SLVI}(F(\bar{x}, \bar{y}), C).$$

*Proof*. This corollary follows from Theorem 5.3 immediately. □

Recall that $K$ is the positive cone in $\mathbb{R}_l^n$; $P_K: \mathbb{R}_l^n \to K$ is the metric projection operator and $\preccurlyeq_K$ is the partial order on $\mathbb{R}_l^n$ induced by $K$. By Definition (3.14), $\preccurlyeq_K$ is also a partial order on $\mathbb{R}_M^n$ induced by $K$. By $P_K(\theta) = \theta$, we consider a special locally variational inequality problem LVI($P_K(\theta, \theta)$, $\mathbb{R}_l^n$). Actually, it is indeed a special global variational inequality problem. The solutions to the problem LVI($P_K(\theta, \theta)$, $\mathbb{R}_l^n$) have the following properties.

**Proposition 5.5**. *Let $P_K: \mathbb{R}_l^n \to K$ be the metric projection operator defined in* (3.14). *For arbitrarily given $y^* \in K$, we have*

$$(x^*, y^*) \in \text{SLVI}(P_K(\theta, \theta), \mathbb{R}_l^n), \text{ for any } x^* \in \mathbb{R}_M^n \text{ with } \theta \preccurlyeq_K x^* \preccurlyeq_K y^*. \tag{5.4}$$

*Consequently,*

$$x^* \in \widehat{D}^*P_K(\theta, \theta)(y^*), \text{ for any } x^* \in \mathbb{R}_M^n \text{ with } \theta \preccurlyeq_K x^* \preccurlyeq_K y^*. \tag{5.5}$$

*Proof.* Let $y^* = (y_1^*, y_2^*, \ldots, y_n^*,) \in K$. Let $x^* = (x_1^*, x_2^*, \ldots, x_n^*,) \in \mathbb{R}_M^n$. Suppose that $\theta \preccurlyeq_K x^* \preccurlyeq_K y^*$. This is, $0 \leq x_i^* \leq y_i^*$, for each $i$ = 1, 2, …, $n$. For any $u = (u_1, u_1, \ldots, u_n) \in \mathbb{R}_l^n$ with $v = (v_1, v_1, \ldots, v_n) =$

$P_K(u) \in K$, we calculate

$$\langle y^*, P_K(u) - P_K(\theta)\rangle - \langle x^*,\ u - \theta\rangle$$
$$= \langle y^*, P_K(u)\rangle - \langle x^*,\ u\rangle$$
$$= \langle y^*, v\rangle - \langle x^*,\ u\rangle$$
$$= \textstyle\sum_{u_i>0}(y_i^* v_i - x_i^* u_i) + \sum_{u_i\leq 0}(y_i^* v_i - x_i^* u_i)$$
$$= \textstyle\sum_{u_i>0}(y_i^* u_i - x_i^* u_i) + \sum_{u_i\leq 0}(y_i^* 0 - x_i^* u_i)$$
$$= \textstyle\sum_{u_i>0}(y_i^* - x_i^*)u_i - \sum_{u_i\leq 0} x_i^* u_i$$
$$\geq 0 \text{ (by } 0 \leq x_i^* \leq y_i^*, \text{ for each } i = 1, 2, \ldots, n).$$

This implies that, for any given $y^* \in K$,

$$(x^*, y^*) \in \text{SLVI}(P_K(\theta,\theta), \mathbb{R}_l^n), \text{ for any } x^* \in \mathbb{R}_M^n \text{ with } \theta \preccurlyeq_K x^* \preccurlyeq_K y^*.$$

This proves (5.4). Then, (5.5) is proved by Theorem 5.3 and (5.4), which reproved this proposition. □

**Acknowledgments** The author is very grateful to Professors Boris S. Mordukhovich and Professor Jen-Chih Yao for their kind communications, valuable suggestions and encouragements in the development stage of this paper.

**Author's declaration:** The author declares that he has no conflict of interest, and the manuscript has no associated data. This work was not supported by any external funding.